%% file: _main_paper.tex
\documentclass[preprint]{elsarticle}

\usepackage{lineno}

\usepackage{amssymb}

\usepackage{blindtext} 
\usepackage{booktabs}
\usepackage{eurosym}
\usepackage{siunitx}
\usepackage{chemformula}
\DeclareSIUnit{\EUR}{\mbox{\euro}}
\DeclareSIUnit{\USD}{\mbox{\$}}
\DeclareSIUnit{\year}{\mbox{y}}
\DeclareSIUnit{\tCO}{\mbox{t\ch{CO2}}}
\DeclareSIUnit{\betaCo}{\mbox{m$^{2\,\beta}$}}
\DeclareSIUnit{\wtpercent}{wt\%}
\DeclareSIUnit{\ct}{ct}
\usepackage{mathrsfs}
\usepackage{float}
\usepackage{threeparttable}
\usepackage{amsmath}

\usepackage{nomencl}
\usepackage{etoolbox}
\usepackage{svg}
\usepackage{graphicx}
\usepackage{float}
\usepackage{comment}
\usepackage{multirow}
\usepackage[justification=centering]{caption}
\usepackage[super]{nth}
\usepackage{url}
\usepackage{soul}
\usepackage{verbatim}

\usepackage{xargs}    

\usepackage{regexpatch}
\makeatletter
\regexpatchcmd\ps@pprintTitle
  {\cE\}\ \cE\}}{\cE\}\cE\}}
  {}{\FailedToPatch}
\makeatother

\begin{document}

\begin{frontmatter}



\title{Unlocking the Potential of Synthetic Fuel Production: Coupled Optimization of Heat Exchanger Network and Operating Parameters of a \mbox{1 MW} Power-to-Liquid Plant}

\journal{Chemical Engineering Science}

\author[inst1]{David Huber\corref{corrauthor}}\ead{david.huber@tuwien.ac.at}
\author[inst1]{Felix Birkelbach}
\author[inst1]{René Hofmann}

\cortext[corrauthor]{Corresponding author}

\address[inst1]{TU Wien, Institute of Energy Systems and Thermodynamics, Getreidemarkt 9/BA, 1060 Vienna, Austria}

\begin{abstract}
\input{00_abstract}
\end{abstract}



\end{frontmatter}


\input{01_introduction}
\input{02_materials}
\input{03_modeling}
\input{04_results}
\input{05_conclusion}

\bibliographystyle{elsarticle-num-names} 
\bibliography{references}

\input{06_statements.tex}

\appendix
\input{A_appendix}
\input{nomenclature}
\end{document}

%% file: 00_abstract.tex
The use of synthetic fuels is a promising way to reduce emissions significantly. To accelerate cost-effective large-scale synthetic fuel deployment, we optimize a novel \SI{1}{\mega\watt} PtL-plant in terms of PtL-efficiency and fuel production costs. For numerous plants, the available waste heat and temperature level depend on the operating point. To optimize efficiency and costs, the choice of the operating point is included in the heat exchanger network synthesis. All nonlinearities are approximated using piecewise linear models and transferred to MILP. Adapting the epsilon constraint method allows us to solve the multi-criteria problem with uniformly distributed solutions on the Pareto front. We improved the lowest production costs of \mbox{$\SI{1.89}{\EUR\per\kilogram}$} and the highest efficiency of \mbox{$\SI{58.08}{\percent}$} from the conventional design process to \mbox{$\SI{1.83}{\EUR\per\kilogram}$} and \mbox{$\SI{61.33}{\percent}$}. By applying the presented method, climate-neutral synthetic fuels can be promoted and emissions can be reduced in the long term.

%% file: 01_introduction.tex
\section{Introduction}
\label{sec:introduction}

The transport sector contributed \mbox{$\SI{7.98}{\giga\tonne}$} of \ch{CO2}, equivalent to about \mbox{$\SI{23}{\percent}$} of global \ch{CO2} emissions, in 2022 \cite{international_energy_agency_co2_2022}. Reducing these emissions is crucial for achieving climate goals. Electrification of the transport sector will enable a significant reduction in climate-damaging emissions. Projections by the IEA indicate that the share of electric vehicles will increase from \mbox{$\SI{1.5}{\percent}$} in 2021 to about \mbox{$\SI{30}{\percent}$} by 2030 and more than \mbox{$\SI{60}{\percent}$} by 2050 \cite{international_energy_agency_global_2022}. Electrification is not possible for the entire transport sector. There is still no viable alternative to liquid fuels for maritime and aviation. Heavy machines such as snow groomers or agricultural machinery are similarly affected by a lack of climate-neutral powertrain systems. Synthetic fuels, or e-fuels, can be used directly in existing combustion engines as a drop-in solution to achieve climate neutrality. Compared to other climate-neutral fuels like hydrogen, using e-fuels offers the advantages of cheap transport and storage costs and existing infrastructure can be continued to be used.

\subsection{Synthetic Fuel Production}
The production of synthetic fuels involves multiple chemical conversion steps of mainly \ch{CO2}, \ch{H2O} and renewable electricity sources like wind, solar or photovoltaic. The \ch{CO2} can be provided through different routes like capturing \ch{CO2} from the atmosphere with direct air capture (DAC), water, biomass, or flue gas \cite{hanggi_review_2019}. Simplified, electrolysis reduces \ch{H2O} to a hydrogen-rich synthesis gas. In a subsequent step, the gas is synthesized with purified \ch{CO2} in a Fischer-Tropsch (FT) reactor. The hydrocarbons from the FT-reactor are then separated into different fractions such as naphtha, diesel and waxes.

Even though the key components such as electrolysis and FT-reactor are already commercially available, only a few plants are still in operation \cite{pinsky_comparative_2020, marchese_energy_2020}. The first commercial PtL-plant build was the Haru Oni in Punta Arenas, Chile \cite{siemens_energy_haru_2023}. Methanol is produced using DAC and a polymer electrolyte membrane (PEM) electrolysis system powered by renewable electricity from a \mbox{$\SI{3.4}{\mega\watt}$} wind turbine. During the pilot phase in December 2022, the first barrel of synthetic fuel could already be filled. The production is to be increased to \mbox{$\SI{130000}{\liter}$} in March 2023. The Norsk e-fuel plant in Mosj\o{}en Norway is expected to produce up to $50$ million liters of kerosene for aviation with renewable electricity from wind and hydropower as early in 2026 \cite{norsk_e-fuel_driving_2023}. The \ch{CO2} required for the solid oxide electrolysis cell (SOEC) is supplied by DAC. According to Norsk e-fuel, three production facilities are expected to produce more than $250$ million liters annually by the end of 2030. The George Olah plant in Svartsengi, Iceland, can produce up to $4000$ tons (about $5$ million liters) of methanol annually \cite{marlin_process_2018}. A nearby geothermal power plant's off-gas provides the \ch{CO2} feedstock. The alkaline water electrolysis is fed entirely renewably from the Icelandic power grid. INERATEC's PtL-plant in Karlsruhe, Germany, is expected to produce about $3500$ tons (about $\SI{4.6}{}$ million liters) of synthetic kerosene and diesel per year, starting in 2023 \cite{ineratec_industrial_2022}. Feedstocks are up to $\SI{10000}{}$ tons of biogenic \ch{CO2} per year and renewable electricity.

The production of e-fuels is still in its earliest stage and significant challenges still need to be overcome \cite{zhao_how_2022}. The central problem is that production costs are currently too high, which means that the fuel cannot be used in an economically viable way for end users \cite{ngando_co_2022}. Ueckerdt et al. \cite{ueckerdt_potential_2021} have estimated production costs of \SI{3,2}{\EUR\per\liter} for gasoline from 2020 to 2050. Ram et al. \cite{ram_powerfuels_2020} have estimated production costs of \SI{1,14}{\EUR\per\liter} for the year 2050. This forecast refers to all FT-products. However, depending on the FT-reactor, the share of gasoline is only about a third, which again leads to an estimate similar to that of Ueckerdt et al. Based on the investment costs of the demonstration plant in Haru Oni, Chile, a current production price of about \SI{50}{\EUR\per\liter} results. In contrast, the average wholesale price of fossil gasoline in 2021 was about \SI{0,5}{\EUR\per\liter} \cite{us_energy_information_administration_us_2023}. A \ch{CO2} price of about \SI{1000}{\EUR\per\tCO} would have to be introduced to create cost equality \cite{ueckerdt_e-fuels_2023}. According to Ueckerdt et al. \cite{ueckerdt_e-fuels_2023}, production costs of about \SI{1}{\EUR\per\liter} will occur in the long term, which will make the fuel interesting for end users. Large and, above all, plants with low production costs and high efficiency must be built to achieve this price target.

\subsection{Plant Design \& Operation}
Complex plants like a PtL-plant are composed of several units. At the design stage of a plant, the size of the units is primarily determined by the intended capacity of the plant. Standard manufacturers' sizes are used to reduce costs by buying off-the-shelf units. To guarantee the proper operation of the plant and the units, process engineers use commercial software such as Aspen HYSYS or IPSEpro. Due to the large number of units interacting with each other, many independent operating parameters result, which makes it challenging to find an optimal operating point empirically. Therefore, the plant’s performance is highly dependent on the knowledge and experience of the process engineers. From an engineering point of view, it is necessary to use optimization methods to utilize the system’s full potential. In contrast to commercial software, operational optimization enables the minimization of certain criteria like costs or emissions to find an optimal operation point. Al-Rashed \& Afrand \cite{al-rashed_multi-criteria_2021}, for example, used a genetic algorithm (GA) to optimize the exergetic and economic efficiency of a combined gas turbine and supercritical \ch{CO2} cycle for power production. They were able to find an optimal operating point with \SI{25.3}{\percent} higher exergetic efficiency and achieved \SI{24.6}{\percent} cost savings resulting from optimized inlet cooling of the compressor. Cao et al. \cite{cao_development_2022} compared biomass gasification and digestion for a combined biomass to power and hydrogen plant. The units’ characteristics were modeled to represent emissions, efficiency and levelized cost of product (LCOP) criteria and optimized with a GA. The results do not indicate a preferred system design. However, modeling the components depending on the operating point could create a holistic basis for decision-making. Hai et al. \cite{hai_-neural_2023} optimized the operation of a solar-geothermal energy system providing electricity and hydrogen using GA. Optimizing the operating point of the components increased the power output by \SI{500}{\kilo\watt} to \SI{4.099}{\mega\watt}. At the same time, the production of \ch{H2} was increased from \SI{8}{\gram\per\second} to \SI{29}{\gram\per\second}. Wang et al. \cite{wang_rolling_2023} optimized the operation strategy of a solar tower power plant using particle swarm optimization. They increased the daily power generation by \SI{13.4}{\percent}. 

The literature cited earlier shows the great potential of optimizing operating characteristics. However, these methods are limited because the heat exchanger network (HEN) design is neglected. The HEN design is a crucial factor for cost and energy savings, especially in plants where much energy is necessary for heating and cooling process streams. Commercial software can be used to simulate an existing HEN. However, the optimal interconnection must be specified by the process engineers. The empirical rules of pinch analysis are often used \cite{linnhoff_introduction_1998}. A promising approach for finding an optimal HEN is the heat exchanger network synthesis (HENS). With HENS, an optimal heat exchanger configuration can be found using mathematical programming subjected to given stream parameters such as temperatures and heat capacity flows. Some representative use cases and the achievable potential for cost and energy savings can be found in the publications of, for example, Yee \& Grossmann \cite{yee_simultaneous_1990}, Escobar et a. \cite{escobar_optimal_2013} and Liu et al. \cite{liu_extended_2022}. 

A significant limitation of HENS is that a defined operating point and, as a result, constant stream parameters are assumed. In cases where the operating point affects the temperatures and flow capacities of the streams, the operating behavior cannot be covered by classical HENS. In contrast, when the operating point is optimized, a predefined configuration of the HEN is necessary. Since the operating point influences the HEN and vice versa, optimizing HEN and the operating point simultaneously will lead to better results. In contrast to flexible HEN design methods, in this paper, we do not optimize the HEN for different operating points. With our method, the operating point and the HEN are optimized simultaneously.

\subsection{Novelty \& Contribution}
In this work, we close a research gap by coupling the optimization of the HEN and operating parameters. Our method allows us to overcome design and operational optimization weaknesses and enforce both approaches' strengths. To exploit the potential of modern MILP solvers, we use piecewise linear approximations and efficient logarithmic coding. Applying an adapted epsilon constraint method enables the consideration of the trade-off between high process efficiency and low production costs. Therefore, decision-makers and process engineers can be provided with a holistic foundation for the process design. 

Our method is demonstrated on a novel \SI{1}{\mega\watt} PtL-plant. As a step towards large-scale industrial production, the goal is to design a medium-sized PtL-plant successfully. For a given electrolyzer size of \SI{1}{\mega\watt}, the aim is to produce as much synthetic fuel as possible cost-effectively and efficiently simultaneously. Therefore, we present how to formulate the efficiency and production cost objective functions. To evaluate our method's potential, the results from coupled optimization are compared to results from the traditional plant design approach, where only the HEN is optimized. The results show that synthetic fuels can be produced with minimal production costs of \mbox{$\SI{1.83}{\EUR\per\kilogram}$} or at the highest efficiency of \mbox{$\SI{61.84}{\percent}$}. The decarbonization of the transport sector can thus be accelerated and emissions can be avoided.

\subsection{Paper Organization}
The novel \SI{1}{\mega\watt} PtL-plant studied is described in Section \ref{sec:systemDescription}. In Section \ref{sec:methods}, we show the basic concept of the adapted HENS superstructure, the piecewise linear approximation, the logarithmic coding approach for an efficient transfer to MILP, and our adaptions to the epsilon constraint method. Since implementing the operating point dependent stream parameters is specifically tailored to the use case, the modeling is discussed in detail in Section \ref{sec:modeling}. In Section \ref{sec:reults}, we show the results of our method in comparison to an optimization of the HEN only.

%% file: 02_materials.tex
\section{Materials \& Methods}
\label{sec:MaterialsMethods}

\subsection{System Description}
\label{sec:systemDescription}
\begin{figure}[H]
    \centering
    \includegraphics[width=1\textwidth]{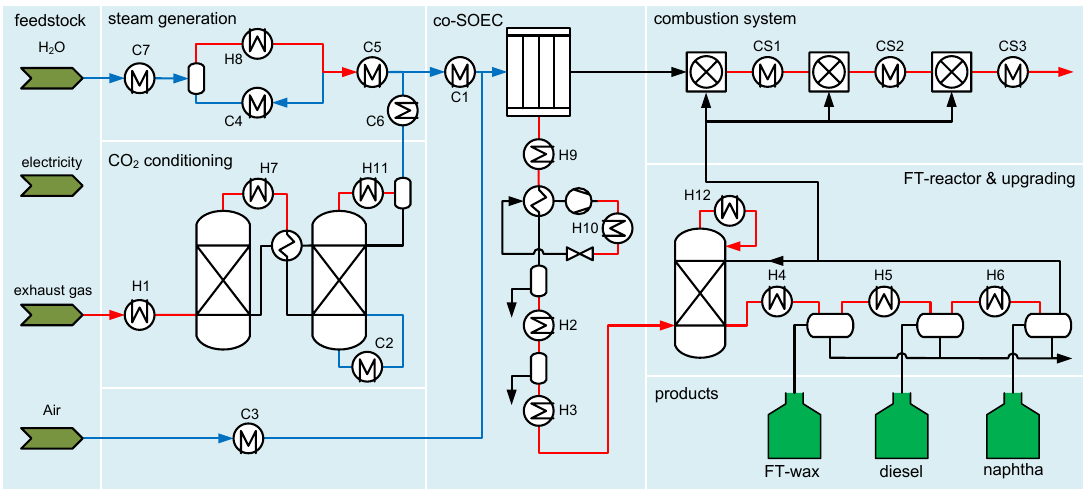}
    \caption{Schematic representation of the \mbox{$\SI{1}{\mega\watt}$} PtL-plant with the five main components: steam generation, \ch{CO2} conditioning, co-SOEC, FT-reactor with upgrading and combustion system.} \label{fig:schematicPfd}
\end{figure}

As part of the \textit{IFE} (de.: Innovation Flüssige Energie, eng.: Innovation Liquid Energy) research project, a novel PtL-plant is designed. The PtL-plant will be designed to a maximum electrolysis capacity of approximately \mbox{\SI{1}{\mega\watt}} and will produce climate-neutral fuels for the transport sector. Water, renewable electricity, exhaust gas from a cement plant and air are used as feedstock. In Figure \ref{fig:schematicPfd}, the PtL-plant is schematically shown with its five main components: steam generation, \ch{CO2} conditioning, high-temperature solid oxide co-electrolysis (co-SOEC), Fischer-Tropsch (FT) reactor with upgrading and combustion system. For simplicity, no valves, pumps, or compressors are shown. Heat exchangers indicate the heat transfer points of the HEN. The cold process streams, which must be heated, are shown in blue. Hot process streams must be cooled down and shown in red. In the following sections, the central components are briefly discussed.

\subsubsection{Steam Generation}
The steam generator is fed with pure water at \SI{20}{\celsius}. The conditioned water is preheated, evaporated and superheated by the heat exchangers in streams C7, C4 and C5.

\subsubsection{\ch{CO2} Conditioning}
Due to the high energy consumption for the supply of \ch{CO2} in sufficient purity, the quality of the \ch{CO2} source is an essential parameter for the design. In this case, exhaust gas from a cement plant is used due to the high \ch{CO2} concentration. The exhaust gas enters the conditioning with \SI{15.17}{\wtpercent} \ch{CO2}, \SI{81.11}{\wtpercent} \ch{N2} and \SI{3.72}{\wtpercent} \ch{H2O} and a temperature of \SI{40}{\celsius}. To control the temperature for adsorption and desorption, the hot streams H7 and H11 must be cooled and the cold stream C2 heated. At the end of the conditioning process, the \ch{CO2} has a purity of \SI{98.73}{\wtpercent} and a low residual content of water.

\subsubsection{Co-SOEC}
The central element of the PtL-process is the co-SOEC. The conditioned \ch{CO2} is superheated within stream C6 and mixed with the superheated steam. Within stream C1, the mixture is further superheated to the operating temperature of the co-SOEC. Both are reformed with preheated air at temperatures between $800$ and $\SI{900}{\celsius}$ to an \ch{H2}-rich gas and \ch{CO}. Before entering the FT reactor, the synthesis gas leaving the co-SOEC is cooled in four stages and condensed constituents are separated. A compressor-based cooling system supports the second cooling stage. 

With a maximum power of approx. \SI{1}{\mega\watt}, the co-SOEC is the largest electricity consumer in the system. The cell voltage has a strong influence on the overall power consumption and subsequently on the efficiency and production costs.

\subsubsection{FT-Reactor \& Upgrading}
In the FT-reactor, the synthetic gas from the co-SOEC is passed over a catalyst at high temperature and pressure. In the subsequent upgrading process, the FT-syncrude is separated into the fractions FT-wax, diesel and naphtha and prepared for final use. The unreacted synthesis gas is partially recirculated and fed to the combustion system. The product properties downstream of the upgrading are given in Table \ref{tab:productProp}.
\begin{table}[H]
    \centering
    \caption{Chemical and physical product properties at $\SI{40}{\celsius}$ and $\SI{101324.97}{\pascal}$ downstream the upgrading.}
    \label{tab:productProp}
    \begin{tabular}{ccrrr}
        \toprule
        $v$ & product & \multicolumn{1}{c}{$h_{\mathrm{prod}}$ / $\SI{}{\mega\joule\per\kilogram}$} & \multicolumn{1}{c}{$\rho_{\mathrm{prod}}$ / $\SI{}{\kilogram\per\cubic\meter}$} & \multicolumn{1}{c}{$\mu_{\mathrm{prod}}$ / $\SI{}{\milli\pascal\per\second}$} \\
        \midrule
        $1$ & FT-wax & $43.887$ & $797.73$ & $6.7477$ \\
        $2$ & diesel & $44.345$ & $748.81$ & $1.5983$ \\
        $3$ & naphtha & $44.676$ & $516.17$ & $0.5893$ \\
        \bottomrule
    \end{tabular}
\end{table}

\subsubsection{Combustion System}
The combustion system serves as an internal hot utility and provides energy to heat the cold process streams. The CS consists of three serially connected combustion chambers. The offgas from the separation, which can no longer be recirculated, is used as fuel. Combustion takes place at high air surplus ($\lambda > 150$ for 1\textsuperscript{st} CS; $\lambda > 30$ for 3\textsuperscript{rd} CS). The first combustion chamber is supplied with the entire exhaust air from the co-SOEC. The other two combustion chambers are each fed with the exhaust gas from the combustion chamber in advance.

Since the CS is used as an internal hot utility, the optimal design is a decisive factor for the efficiency of the HEN and the overall process. The selection of the combustion parameters inlet and outlet temperature into the combustion chamber and fuel mass flow are, in contrast to other processes, not closely linked to technical limitations. Only the material-specific temperature limit of \SI{900}{\celsius} must not be exceeded. The inlet temperature of the air mainly influences the outlet temperature from the combustion chamber and can be adjusted by the fuel quantity. The inlet temperature results from the streams connected within the HEN. Due to the many freely selectable parameters and the strong influence on the efficiency of the overall process, the design of the CS is a significant challenge for process engineers.

\subsection{Methods}
\label{sec:methods}
\begin{figure}[H]
    \centering
    \includegraphics[width=1\textwidth]{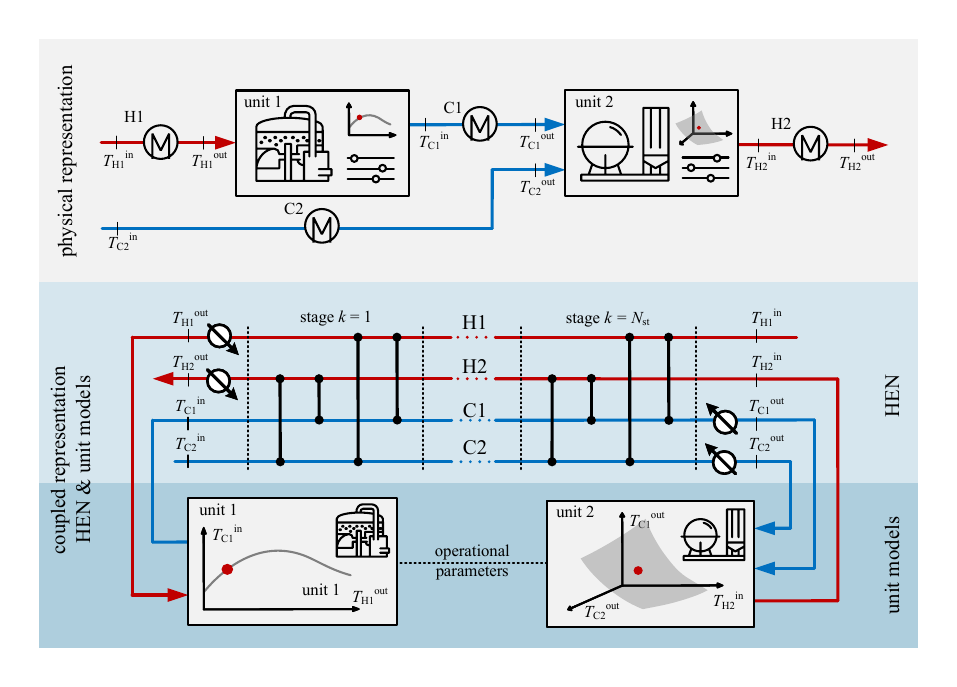}
    \caption{Schematic representation of the method for two units with two hot and two cold streams.}
    \label{fig:schematic}
\end{figure}

In this section, we present the underlying idea behind the coupled optimization of the HEN and operating point. The concept is shown schematically in Figure \ref{fig:schematic}. The physical system composed of two units with two hot and two cold streams is shown at the very top of Figure \ref{fig:schematic}. For simplicity, the heat exchanger network is not illustrated. The operating characteristics of each unit can be reduced to multiple independent operational parameters. It is assumed that a change of the operating parameters leads only to a change of the stream parameters $T^{\textrm{in}}$, $T^{\textrm{out}}$ and $F$. Depending on the process, all, none, or only some of the three stream parameters can be affected.

In the middle of Figure \ref{fig:schematic}, the HEN with $N_{\textrm{st}}$ stages and all possible interconnections are shown for the same system. This graph-theoretic representation of the HEN is based on the superstructure formulation of Yee \& Grossman \cite{yee_simultaneous_1990}. The heat exchange between the hot and cold streams can occur in $N_{\textrm{st}}$ stages with stream splits. The hot and cold utilities are located at the stream ends. 

The coupling between the optimization of the HEN design and the operating point is achieved by implementing models that describe the operating characteristics of the units. In Figure \ref{fig:schematic}, at the bottom, the models of the two units are shown schematically. Each unit has an independent operating variable; in this case, only inlet or outlet temperatures change  formulated with variables. The model of unit 1 represents the correlation between the outlet temperature of stream H1 and the inlet temperature of stream C1. The model of unit 2 illustrates the correlation between the outlet temperatures of stream C1 and C2 along with the inlet temperature of stream H2.  The objective links the operating variables. The coupled HEN design problem with variable stream parameters is solved using the adapted superstructure formulation of Huber et al. \cite{huber_hens_2023}.

\subsubsection{Linearization}
\label{sec:linearization}
All nonlinearities of the HENS are piecewise linear approximated with superpositioned planes in the two-dimensional. Plane Simplices are used for three-dimensional correlations. Detailed information regarding the methodology can be found in the paper by Huber et al. \cite{huber_hens_2023}.

Analogously, the units' operating and stream parameters correlations are linearly approximated. The advantage of this method is that the characteristics of the units can be represented with the help of black-, grey- or white-box models. The resulting flexibility in choosing the modeling approach reduces limitations regarding data availability and increases the proposed method's applicability.

\subsubsection{Transfer to MILP}
\label{sec:transferToMILP}
All piecewise linear functions are transferred to MILP with the least possible number of binary variables to reduce computation time. One-dimensional, mainly convex curved functions are transferred to MILP without binary variables. All other functions require the use of binary variables. Applying a logarithmic coding approach, according to Vielma and Nemhauser \cite{vielma_modeling_2011}, can reduce the number of binary variables to a minimum. Further information about the transfer to MILP can be taken from \cite{huber_hens_2023}.

\subsubsection{Multi-Objective Optimization}
\label{sec:MOO}

The two-objective optimization problem is solved with an adapted epsilon constraint method. As shown in Equation \eqref{eq:epsilon}, one objective function is minimized and the second is constrained with an upper and lower bound. The epsilon parameters are chosen to be in the range of $f_{2}^{\textrm{min}} \le \varepsilon_i \le f_{2}^{\textrm{max}}$.
\begin{equation}
\label{eq:epsilon}
\begin{split}
    &\textrm{min} \,\,\, f_{1}(\boldsymbol{x}) \\
    &\textrm{s.t.} \,\,\, \varepsilon_{i} \le f_{2}(\boldsymbol{x}) \le \varepsilon_{i+1} \quad i=1,\dots,m
\end{split}
\end{equation}
In contrast to the conventional epsilon constraint method, overhanging regions of the Pareto front can be covered. Equidistantly distributed points on the Pareto front can be calculated with minimal computational effort.

%% file: 03_modeling.tex
\section{Modeling} 
\label{sec:modeling}
The models are created based on data from steady-state process simulations by project partners using Aspen HYSIS. Modeling is based on the following assumptions:
\begin{itemize}
    \item The system was simulated with seven different cell voltages between lower ($U_{\textrm{cell}}^{\textrm{min}} = \SI{1.275}{\volt}$) and the upper technical limit ($U_{\textrm{cell}}^{\textrm{max}} = \SI{1.305}{\volt}$).
    \item The stream parameters $T^{\textrm{in}}$, $T^{\textrm{out}}$ and $F$ are only dependent on the cell voltage \mbox{$U_{\textrm{cell}}$}.
    \item The stream parameters of the CS are independent of the cell voltage. Temperatures and heat capacity flows are limited by the amount of offgas available.
    \item The sizes and parameterization of the units are independent of the cell voltage. Identical system costs are assumed.
\end{itemize}

\subsection{Feedstock}
\label{sec:feedstockModels}
\begin{figure}[H]
    \centering
    \includegraphics[width=1\textwidth]{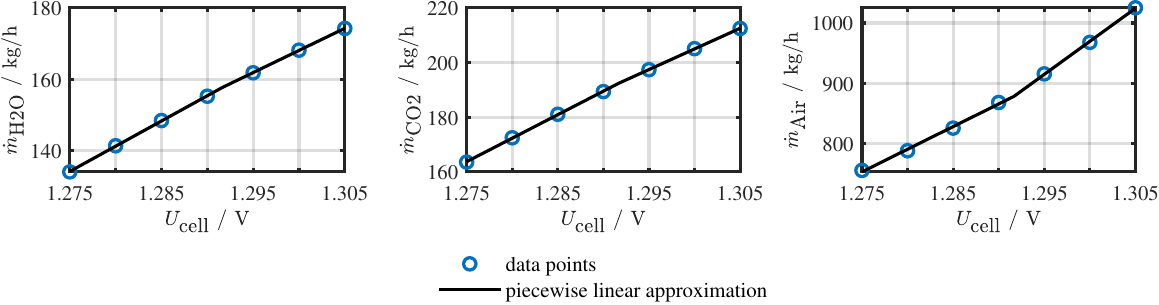}
    \caption{Piecewise-linear approximations for feedstocks as a function of the cell voltage $U_{\textrm{cell}}$. Left: 2 lines, \mbox{$RMSE = \SI{0.25}{\percent}$}. Center: 2 lines, \mbox{$RMSE = \SI{0.25}{\percent}$}. Right: 2 lines, \mbox{$RMSE = \SI{0.75}{\percent}$}.}
    \label{fig:feedstock}
\end{figure}

In addition to climate-neutral electricity, the primary feedstocks of the PtL-process are \ch{H2O}, \ch{CO2} and air. All three mass flows depend only on the cell voltage and implemented with variables. Figure \ref{fig:feedstock} shows the piecewise linear models of the feedstock flows.

\subsection{System Power \& FT-products}
\label{sec:perfModels}
\begin{figure}[H]
    \centering
    \includegraphics[width=1\textwidth]{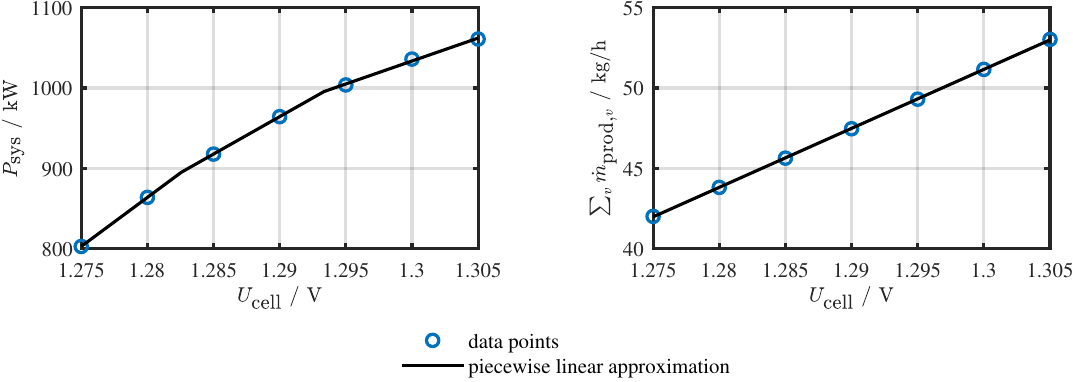}
    \caption{Piecewise-linear approximations of $P_{\textrm{sys}}$ (left) and massflow of all FT-products (right) as a function of the cell voltage $U_{\textrm{cell}}$. Left: 3 lines, \mbox{$RMSE = \SI{0.43}{\percent}$}. Right: 1 line, \mbox{$RMSE = \SI{0.19}{\percent}$}}
    \label{fig:perf}
\end{figure}

The parameters which significantly influence the overall process's performance are the system power $P_{\textrm{sys}}$ and the massflow of FT-products \mbox{$\sum_v \dot m_{\textrm{prod,}v}$}. Figure \ref{fig:perf} shows the simulated data points and the model with piecewise linear approximated lines as a function of cell voltage $U_{\textrm{cell}}$. Both $P_{\textrm{sys}}$ and product output increase with higher cell voltage. The model on the left side in Figure \ref{fig:perf} describes the system’s power consumption. $P_{\textrm{sys}}$ represents the required power to run the co-SOEC, circulation pumps, valves, and control equipment, and to cover losses. The power consumption of utilities is considered in the modeling of the objectives in Section \ref{sec:objectives}. Figure \ref{fig:perf} on the right shows the product flow downstream of the FT-reactor as the sum of the fractions FT-wax, diesel and naphtha.

\begin{figure}[H]
    \centering
    \includegraphics[width=1\textwidth]{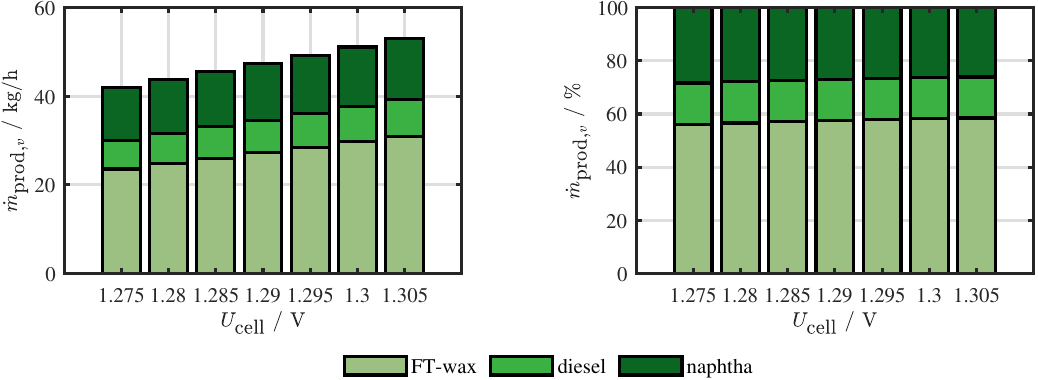}
    \caption{Absolute and relative product composition as a function of cell voltage $U_{\textrm{cell}}$.}
    \label{fig:product}
\end{figure}

The composition of the fractions changes depending on the cell voltage. Figure \ref{fig:product} shows the absolute and relative share of the product flows. The percentage of FT-wax increases with increasing cell voltage. With an almost unchanged share of diesel, the share of naphtha is reduced simultaneously. Since we subsequently relate the production costs to all FT-products, the relative share of the fractions is not considered further.

\subsection{Streams}
\label{sec:streamModels}
\begin{figure}[H]
    \centering
    \includegraphics[width=1\textwidth]{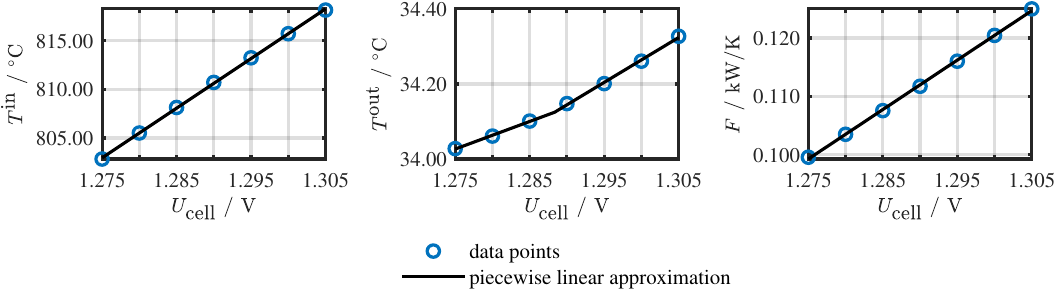}
    \caption{Piecewise-linear approximation of Stream H9. Left: 1 segment, \mbox{$RMSE = \SI{0.55}{\percent}$}. Center: 2 segments, \mbox{$RMSE = \SI{0.90}{\percent}$}. Right: 1 segment, \mbox{$RMSE = \SI{0.83}{\percent}$}.}
    \label{fig:H9fig}
\end{figure}

The stream parameters $T^{\textrm{in}}$, $T^{\textrm{out}}$ and $F$ are dependent on the cell voltage. For each parameter, a piecewise linear approximation is generated. Figure \ref{fig:H9fig} shows the results from process simulation and the piecewise linear approximation for stream H9. In this case, all parameters increase with increasing cell voltage. All other stream models can be found in the supplementary material.

Table \ref{tab:streamData} summarizes the stream parameters’ bounds for all streams. The stream parameters with brackets are variable and implemented depending on the cell voltage. Values without brackets are independent of cell voltage and constant. Especially for processes that require defined inlet and outlet conditions, the temperatures remain constant and only the flow capacity changes. This can be seen in the stream data of streams H12 and H4-H6 of the FT-reactor and upgrading. Stream parameters in square brackets are variable and depend on the cell voltage. The minimum and maximum values from the process simulation are given. The overall heat transfer coefficients were assumed to be \mbox{$U = \SI{0.5}{\kilo\watt\per\meter\squared\per\kelvin}$} for all streams.
\begin{table}[H]
    \centering
    \caption{Stream data: inlet, outlet temperature and flow capacity. The bounds for variable defined stream parameters are shown in square brackets.}
    \label{tab:streamData}
    \begin{tabular}{crrr}
        \toprule
        Stream & $T^{\mathrm{in}}$ / $\SI{}{\celsius}$ & $T^{\mathrm{out}}$ / $\SI{}{\celsius}$ & $F$ / $\SI{}{\kilo\watt\per\kelvin}$ \\
        \midrule 
        H1 & $40.0$ & $35.0$ & $[1.71, 2.16]$\\
        H2 & $[127.9, 131.1]$ & $[34.0, 35.0]$ & $[0.09, 0.12]$\\
        H3 & $[169.8, 174.1]$ & $[34.0, 35.0]$ & $[0.09, 0.12]$\\
        H4 & $210.0$ & $190.0$ & $[0.27, 0.28]$\\
        H5 & $190.0$ & $120.0$ & $[0.56, 0.58]$\\
        H6 & $120.0$ & $30.0$ & $[0.48, 0.50]$\\
        H7 & $[45.4, 57.0]$ & $31.0$ & $[2.35, 2.95]$\\
        H8 & $138.9$ & $137.9$ & $[59.60, 94.40]$\\
        H9 & $[805.2, 825.5]$ & $[34.0, 35.0]$ & $[0.10, 0.13]$\\
        H10 & $[49.5,50.7]$ & $[34.0, 35.0]$ & $[0.65, 0.88]$\\
        H11 & $101.8$ & $30.0$ & $[0.51 ,0.64]$\\
        H12 & $190.0$ & $188.0$ & $[76.88, 80.45]$\\
        C1 & $[318.0, 319.2]$ & $[825.0, 870.5]$ & $[0.14, 0.18]$\\
        C2 & $116.9$ & $124.2$ & $[20.02, 25.12]$\\
        C3 & $[57.3, 58.8]$ & $825.0$ & $[0.25, 0.33]$\\
        C4 & $137.9$ & $139.9$ & $[105.77, 142.64]$\\
        C5 & $138.9$ & $[426.6, 449.4]$ & $[0.10, 0.11]$\\
        C6 & $35.0$ & $[115.9, 145.4]$ & $[0.05, 0.06]$\\
        C7 & $20.3$ & $[189.5, 199.6]$ & $[0.15, 0.21]$\\
        \bottomrule
    \end{tabular}
\end{table}

\subsection{Combustion System}
The stream parameters of the CS are independent of the cell voltage; thus, no piecewise linear approximation is needed. Table \ref{tab:streamDataOC} summarizes all stream parameters and heat transfer coefficients. The exhaust stream's inlet temperature corresponds to the combustion chamber's outlet temperature. An upper technical limit of $\SI{900}{\celsius}$ was chosen for all three streams. The outlet temperature can vary between $\SI{100}{\celsius}$ and $\SI{890}{\celsius}$. The boundaries of the heat capacity flows were selected to guarantee that a sufficient amount of offgas is available for combustion. All overall heat transfer coefficients were assumed to be \mbox{$U = \SI{0.5}{\kilo\watt\per\meter\squared\per\kelvin}$}.
\begin{table}[H]
    \centering
    \caption{Stream data with limits for inlet, outlet temperature and flow capacity.}
    \label{tab:streamDataOC}
    \begin{tabular}{crrr}
        \toprule
        Stream & $T^{\mathrm{in}}$ / $\SI{}{\celsius}$ & $T^{\mathrm{out}}$ / $\SI{}{\celsius}$ & $F$ / $\SI{}{\kilo\watt\per\kelvin}$\\
        \midrule
        CS1 & $900.0$ & $[100, 890]$ & $[59.60, 94.40]$\\
        CS2 & $900.0$ & $[100, 890]$ & $[0.10, 0.13]$\\
        CS3 & $900.0$ & $[100, 890]$ & $[0.65, 0.88]$ \\
        \bottomrule
    \end{tabular}
\end{table}

\subsection{Objectives}
\label{sec:objectives}
Considering objective functions describing efficiency and costs simultaneously makes sense when optimizing energy conversion systems \cite{al-rashed_multi-criteria_2021}. These objective functions are commonly antagonistic, allowing multiple optimal solutions to be represented as a Pareto front. The fuel production costs are minimized and the PtL-efficiency is maximized. The two objective functions reflect the trade-off between efficiency and productivity.

In Bezug auf Gleichung 1  
\subsubsection{PtL-Efficiency}
\label{sec:objPtL}
The objective to maximize the PtL-efficiency $\eta_{\textrm{PtL}}$ is defined by Equation \eqref{eq:objEta} as the ratio of chemically bounded energy to electrical energy input. 

\begin{equation}
    \max \eta_{\textrm{PtL}} = \frac{\dot H_{\mathrm{prod}}}{P_{\textrm{el}}} =
    \frac{ \sum_{v}^{} \dot{m}_{\mathrm{prod,}v} \, h_{\mathrm{prod,}v}}{
    P_{\textrm{sys}} + \sum_{j} \varepsilon_{\textrm{hu}} \, q_{\mathrm{hu},j} + \sum_{i} \varepsilon_{\textrm{cu}} \, q_{\mathrm{cu},i}}
\label{eq:objEta}
\end{equation}

The numerator describes the chemically bounded energy of the FT-products downstream of the separation. The denominator represents the total electrical energy input $P_{\textrm{el}}$. Both the hot and cold utilities are electrified. The coefficient of performance $\varepsilon$ describes the electrical-to-thermal energy input ratio.

\begin{figure}[H]
    \centering
    \includegraphics[width=1\textwidth]{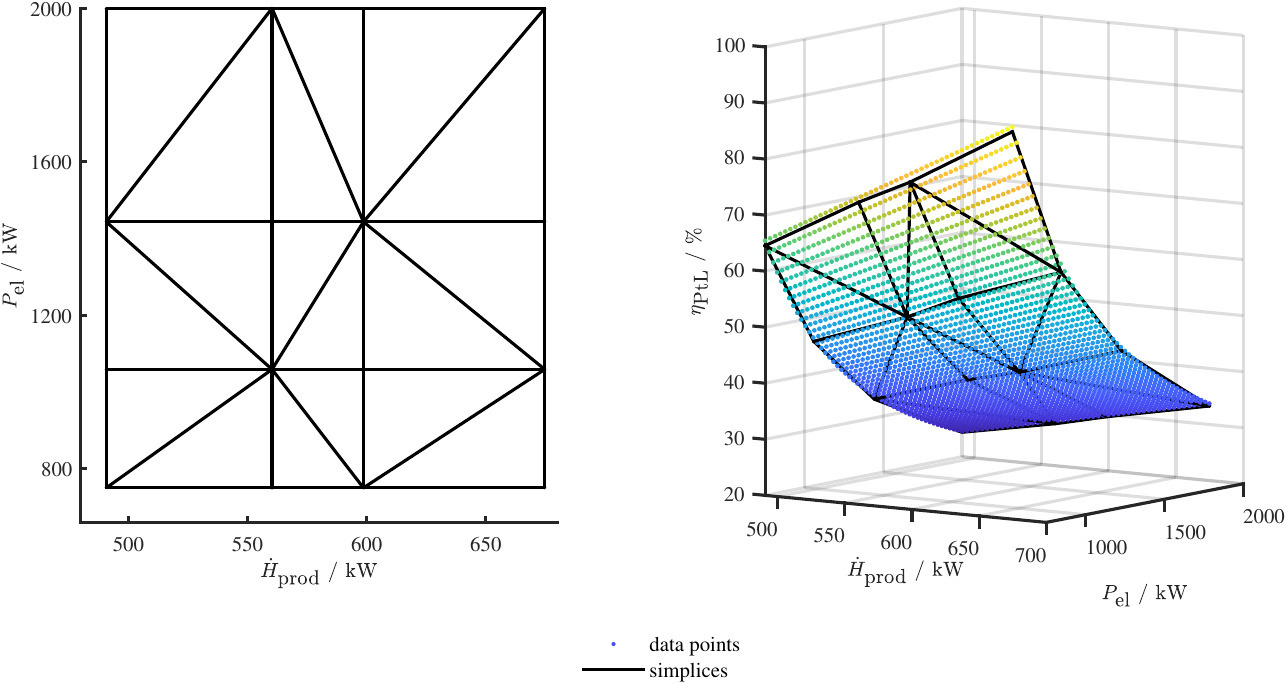}
    \caption{Piecewise-linear approximation with 18 simplices of the PtL-efficiency $\eta_{\textrm{PtL}}$ as a function of chemically bounded energy in the produced fuels $\dot H_{\mathrm{fuel}}$ and electricity demand $P_{\textrm{el}}$. \mbox{$RMSE = \SI{0.61}{\percent}$}.}
    \label{fig:objEtaPtL}
\end{figure}

The Ptl-efficiency from Equation \eqref{eq:objEta} forms a non-linear correlation of the optimization variables for the chemically bounded energy in the product $\dot H_{\mathrm{prod}}$ and the electrical energy required $P_{\textrm{el}}$. The points in Figure \ref{fig:objEtaPtL} on the right show the non-linear function of the PtL-efficiency. Using MILP requires a piecewise-linear approximation. In this case, we used simplices on a regular grid with $4$ x $4$ points, see Figure \ref{fig:objEtaPtL} on the left. With $18$ simplices, the objective can be approximated with sufficient accuracy at an RMSE of $\SI{0.61}{\percent}$.

\subsubsection{Production Costs}
\label{sec:objFuel}

The objective to minimize the specific production costs $c_{\textrm{prod}}$ is defined according to Equation \eqref{eq:objFuel} and describes the ratio of total annual costs to product output. 

\begin{equation}
    \min c_{\textrm{prod}} = \frac{ \mathit{TAC}}{\sum_{v}^{} t \, \dot{m}_{\mathrm{prod,}v}} = \frac{\mathit{CAPEX} + \mathit{OPEX}}{\sum_{v}^{} t \, \dot{m}_{\mathrm{prod,}v}}
\label{eq:objFuel}
\end{equation}

The $TAC$ are composed of annual capital expenses $\mathit{CAPEX}$ and operational expenditures $\mathit{OPEX}$. According to Equation \eqref{eq:CAPEX1}, the $\mathit{CAPEX}$ comprises the investment costs for the system $\mathit{CAPEX}_{\textrm{sys}}$ and the heat exchanger network $\mathit{CAPEX}_{\textrm{HEN}}$. The investment costs for the system are calculated according to Equation \eqref{eq:CAPEX2} and include all relevant investments related to the PtL-plant, excluding the costs for the heat exchanger network. The annualized investment costs for the heat exchanger network are calculated according to Yee \& Grossman \cite{yee_simultaneous_1990} with Equation \eqref{eq:CAPEX3}. Fixed costs for all heat exchangers and variable costs proportional to the heat exchanger area are considered. The investment annualization factor can be expressed as \mbox{$\mathit{AF}_{\mathrm{inv}} = \nicefrac{1}{a}$}. Where $a$ is the depreciation period in years. In contrast to linear depreciation, other options are also possible \cite{Halmschlager2022}.

\begin{equation}
    \mathit{CAPEX} = \mathit{CAPEX}_{\textrm{sys}} + \mathit{CAPEX}_{\textrm{HEN}}
\label{eq:CAPEX1}
\end{equation}

\begin{equation}
    \mathit{CAPEX}_{\textrm{sys}} = \mathit{AF}_{\mathrm{inv}} \, C_{\textrm{sys}}
\label{eq:CAPEX2}
\end{equation}

\begin{equation}
\begin{split}
    \mathit{CAPEX}_{\textrm{HEN}} &= \mathit{AF}_{\mathrm{inv}} \left[ \underbrace{\sum_{i}\sum_{j}\sum_{k}{c_\mathrm{v,hex}\left(\frac{q_{ijk}}{U_{ij}\, \mathit{LMTD}_{ijk}} \right)^{\beta}}}_\text{variable HEX stream costs} \right. \\
    &+ \left.\underbrace{\sum_{i}{c_\mathrm{v,hex} \left(\frac{q_{\mathrm{cu},i}}{U_{\mathrm{cu},i}\, \mathit{LMTD}_{\mathrm{cu},i}} \right)^{\beta}}}_\text{variable HEX cold utility costs}  + \underbrace{\sum_{j}{c_\mathrm{v,hex} \left(\frac{q_{\mathrm{hu},j}}{U_{\mathrm{hu},j}\, \mathit{LMTD}_{\mathrm{hu},j}} \right)^{\beta}}}_\text{variable HEX hot utility costs} \right. \\
    &+ \left. \underbrace{\sum_{i}\sum_{j}\sum_{k}{c_{\mathrm{f,hex}}\,z_{ijk}} + \sum_{i}{c_{\mathrm{f,hex}} \, z_{\mathrm{cu},i}} + \sum_{j}{c_{\mathrm{f,hex}} \, z_{\mathrm{hu},j}}}_\text{fixed investment costs hex} \right]
\label{eq:CAPEX3}
\end{split}
\end{equation}

According to Equation \eqref{eq:OPEX}, the operational expenses are composed of feedstock and electricity costs. The OPEX depend on the annual full load hours $t$ and are usually depreciated within one year, which results in an operational annualization  factor of \mbox{$\mathit{AF}_{\mathrm{op}} = 1$}.

\begin{equation}
\begin{split}
    \mathit{OPEX} &= \mathit{AF}_{\mathrm{op}} \, t \left\{ 
    \underbrace{\sum_{w} c_{\textrm{f}} \, \dot V_{\textrm{f}}}_\text{feedstock costs}  \right. \\
    &+ \left. \underbrace{ c_{\textrm{el}} \, \left[  P_{\textrm{sys}} +  \left( \sum_{j} \varepsilon_{\textrm{uh}} \, q_{\mathrm{uh},j} +  \sum_{i} \varepsilon_{\textrm{uc}} \, q_{\mathrm{uc},i} \right) \right]}_\text{electricity costs}  \right\}
\label{eq:OPEX}
\end{split}
\end{equation}

\begin{figure}[H]
    \centering
    \includegraphics[width=1\textwidth]{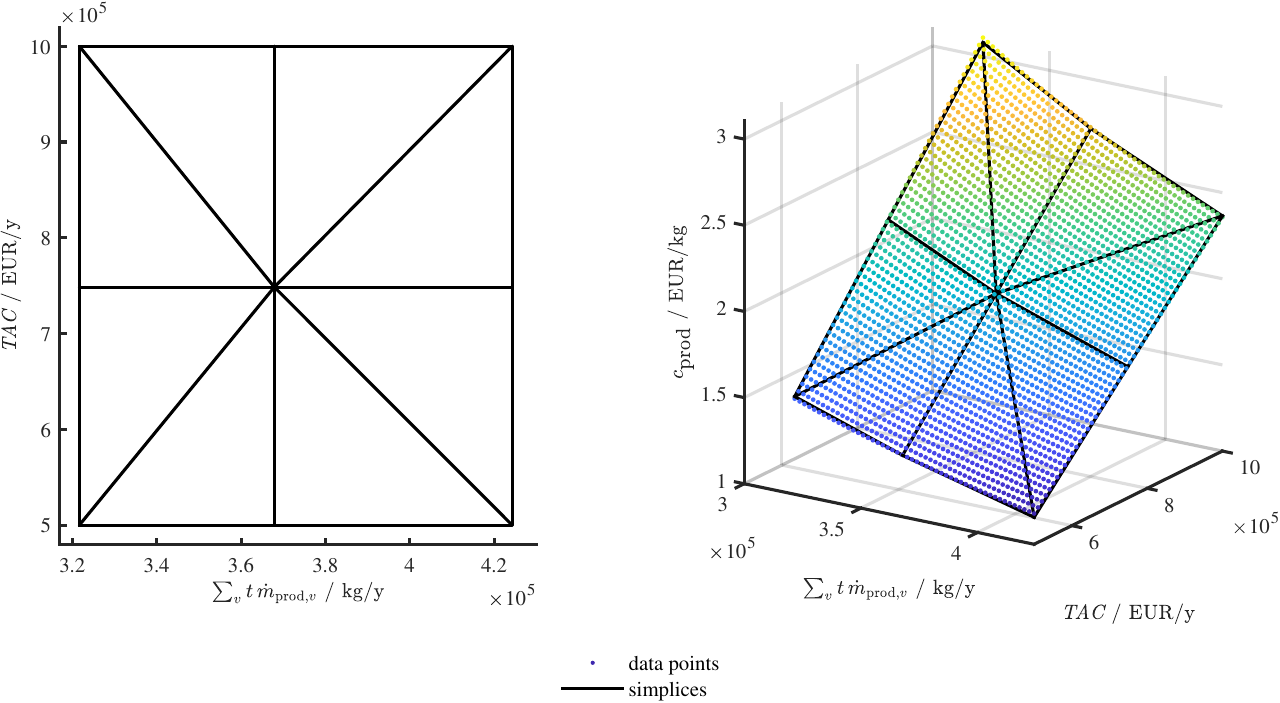}
    \caption{Piecewise-linear approximation with eight simplices of production costs $c_{\textrm{prod}}$ as a function of the total annual costs $TAC$ and that total product flow $\sum_{v}^{} \dot{m}_{\mathrm{prod,}v}$. \mbox{$RMSE = \SI{0.37}{\percent}$}.}
    \label{fig:objFuel}
\end{figure}

The nonlinear objective of the minimum fuel production costs is piecewise-linear approximated and transferred to MILP. Figure \ref{fig:objFuel} shows the approximation in its valid domain on a $3$ x $3$ grid. With only eight simplices, an $\mathit{RMSE}$ of $\SI{0.37}{\percent}$ can be achieved. 

The three nonlinear terms of the variable heat exchanger costs in Equation \eqref{eq:OPEX} are piecewise-linear approximated analogously to the procedure of Huber et al. \cite{huber_hens_2023}. Accordingly, the stream and utility heat exchanger area correlations are approximated by superpositioned planes and transferred to MILP without additional binary variables.

\subsection{Simulation Parameters}
\label{sec:SimParameters}

\subsubsection{Costs}
As the number of full load hours of the plant increases, the net production costs decrease \cite{adelung_global_2022}. Like the authors in \cite{adelung_global_2022,herz_economic_2021,dieterich_power--liquid_2020}, we also assume \mbox{$t = \SI{8000}{\hour\per\year}$} full load hours per year. All components are depreciated linearly within 20 years resulting in an annualization factor of \mbox{$\mathit{AF}_{\mathrm{inv}} = \nicefrac{1}{20}$}. 

The fixed investment costs have been estimated by the project partners responsible for the economic viability to be $C_{\textrm{sys}} = \SI{10000000}{\EUR}$. This estimate assumes that the cost of the central components will decrease significantly due to technological advances \cite{ueckerdt_potential_2021}. In the papers of G. Herz et al. \cite{herz_economic_2021} and D.H. König et al. \cite{konig_techno-economic_2015}, costs in the range of $\SI{9600000}{\EUR}$ and $\SI{22000000}{\EUR}$ have been predicted, which legitimizes the intra-project estimation.

The feedstock for the production of synthetic fuels can be reduced to \ch{H2O}, \ch{CO2} rich exhaust gas and air in a simplified form. Table \ref{tab:feestockCosts} lists the corresponding cost factors $c_{\textrm{f}}$.
\begin{table}[H]
    \centering
    \caption{Feedstock cost factors.}
    \label{tab:feestockCosts}
    \begin{threeparttable}
    \begin{tabular}{cccl}
        \toprule
        $f$ & feedstock & $c_{\textrm{f}} / \SI{}{\EUR\per\tonne}$ & reference \\
        \midrule
        $1$ & \ch{H2O} & $3.54$ & mean for Europe in \cite{eureau_europes_2021} \\
        $2$ & \ch{CO2} & $50^1$ & EU ETS 2026-30 \cite{ieta_ghg_2022}\\
        $3$ & Air & $0^2$ &  \\
        \bottomrule
    \end{tabular}
    \begin{tablenotes}[flushleft]
    \footnotesize
        \item $^1$ \ch{CO2} treatment is considered in $C_{\textrm{sys}}$ and $P_{\textrm{sys}}$.
        \item $^2$ Air is available free of charge. The process was designed to use air with ambient conditions.
    \end{tablenotes}
    \end{threeparttable}
\end{table}

The goal of producing \ch{CO2}-neutral fuels can only be achieved by using electricity produced in a \ch{CO2}-neutral way. The electricity required for this will be generated from a mix of wind and solar PV. At a projected plant location in Europe, an average electricity price of \mbox{$c_{\textrm{el}}=\SI{20}{\EUR\per\mega\watt\per\hour}$} is assumed, according to J.L.L.C.C. Janssen et al. \cite{janssen_country-specific_2022}.

\subsection{Heat Exchanger Network}
The heat exchanger network was defined with \mbox{$N_{\textrm{st}}=3$} stages for heat exchange and a minimum temperature difference of \mbox{$\Delta T_{\textrm{min}}=\SI{1}{\kelvin}$}.

Air coolers are used as cold utilities. Since the environment is used as a heat sink, a coefficient of performance of \mbox{$\varepsilon_{\textrm{uc}} = 0.05$} is assumed to cover evaporation losses and mechanical work. Electrical heating elements with a coefficient of performance of \mbox{$\varepsilon_{\textrm{uh}} = 1.05$} are used as hot utilities.

The cost parameters for the heat exchangers in the heat exchanger network have been determined based on the DACE Price Booklet \cite{noauthor_dace_2021}. We have selected stainless plate heat exchangers made of AISI 316. The cost parameters are $c_\mathrm{f,hex} = \SI{1013.6}{\EUR\per\year}$ and $c_\mathrm{v,hex} = \SI{61.8}{\EUR\per\betaCo\per\year}$ with a degressive cost exponent of $\beta = \SI{0.8}{}$.

\subsection{Piecwise-linear Approximation}
Superpositioned lines, planes or simplices were added to all piecewise-linear models until an RMSE of less than $\SI{1}{\percent}$ was achieved. When approximating with simplices, care was taken to fully utilize the logarithmic encoding with respect to the selected grid points.

\subsection{Implementation}
\label{sec:implementation}
The multiobjective optimization problem was solved using the adapted epsilon constraint method shown in Equation \eqref{eq:epsilon}. For that, $\eta_{\textrm{PtL}}$ was implemented as $-f_{1}(\boldsymbol{x})$ and $c_{\textrm{prod}}$ as $f_{2}(\boldsymbol{x})$.

All optimization problems were coded in MATLAB 2020b \cite{the_mathworks_inc_matlab_2022}. YALMIP R20210331 has been used as an interface between MATLAB and the MILP solver \cite{lofberg_toolbox_2004}. Gurobi 10.0.0 is used as the MILP Solver \cite{gurobi_optimization_llc_gurobi_2023}. A MIP of less than \mbox{$\SI{1}{\percent}$} was defined as a termination criterion for the solver. The calculations were performed on a 64-core server (AMD EPYC 7702P) with $\SI{265}{\giga\byte}$ of RAM.

%% file: 04_results.tex
\section{Results}
\label{sec:reults}
The coupled optimization results are compared with two conventionally designed plants without optimization to quantify the presented method's potential. The following assumptions are made for the conventional plant design:  
\begin{itemize}
    \item The two operating points are set at the extreme cell voltage values.
    \item The HEN is designed empirically following the PINCH rules.
    \item The stream parameters of the CS are chosen empirically.
\end{itemize}

\begin{figure}[H]
    \centering
    \includegraphics[width=1\textwidth]{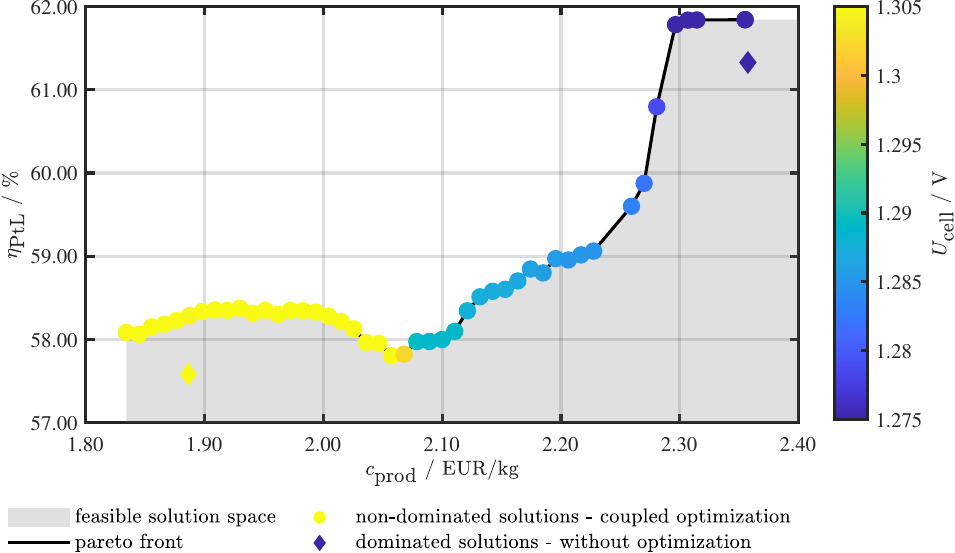}
    \caption{Non-dominated solutions and Pareto front of coupled optimization (circles). Dominated solutions of the conventional approach (diamonds).}
    \label{fig:reults}
\end{figure}

Figure \ref{fig:reults} shows the two solutions of the conventional design without optimization and the Pareto front of the coupled optimization as a function of the cell voltage. The different cell voltages are color-coded. The Pareto front is constructed with 42 points. On average, a solver time of $\SI{465.7}{\second}$ per point was achieved. All $42$ points were calculated within $23.8$ hours. More details can be found in \ref{sec:solvertime}. The gaps in areas of high production costs result from the solver timeout. Within $12$ hours, no solutions with a MIP gap of less than $\SI{1}{\percent}$ could be found. A high cell voltage generally leads to low efficiency and production costs. A low cell voltage, on the other hand, leads to increased efficiency and, at the same time, high production costs. The two objective functions are intensely antagonistic, at least near the optimum.

Table \ref{tab:results} shows the characteristic values for the empirical design approach without optimization. Further, the characteristic values at the corners of the Pareto front from the coupled optimization are shown. The associated stream plots can be found in \ref{sec:streamplots}.

Production costs at a cell voltage of $\SI{1.275}{\volt}$ are almost identical at \mbox{$\SI{2.355}{\EUR\per\kilogram}$} compared to \mbox{$\SI{2,358}{\EUR\per\kilogram}$}. The efficiency can be increased by 0.514 percentage points ($\SI{0.84}{\percent}$). Both HENs use $27$ heat exchangers. Whereas in the empirical design, one is used as hot utility HEX. Further, it is noticeable that the major improvement results from the lower power of the cold utilities. The power can be reduced from $\SI{222.25}{\kilo\watt}$ to $\SI{155.62}{\kilo\watt}$.

At a cell voltage of $\SI{1.305}{\volt}$, the production cost can be reduced by $\SI{5.22}{\ct}$ ($\SI{2.77}{\percent}$) while increasing the efficiency by $0.504$ percentage points ($\SI{0.88}{\percent}$). Compared to the HEN of the coupled optimization with $20$ HEX, the empirical design requires $34$ HEX. The lower cost of HEN is a crucial driver for the cost savings. It should also be pointed out that only one of the three CS is required for the optimized HEN; see Figure \ref{fig:streamplotopt13005} in \ref{sec:streamplots}. This may result in cost reductions in system costs. However, these will not be discussed further in this paper.

\begin{table}[H]
    \centering
    \caption{Comparison of extreme values from coupled optimization and empirical design without optimization.}
    \label{tab:results}
    \begin{tabular}{c c r r r r}
        \toprule
         variable & unit &  \multicolumn{1}{c}{empirical} & \multicolumn{1}{c}{optimization} & \multicolumn{1}{c}{empirical} & \multicolumn{1}{c}{optimization} \\
        \midrule
        $U_{\textrm{cell}}$ & $\SI{}{\volt}$ & $\SI{1.275}{}$ & $\SI{1.275}{}$ & $\SI{1.305}{}$ & $\SI{1.305}{}$ \\
        $c_{\textrm{prod}}$ & $\SI{}{\EUR\per\kilogram}$ & $\SI{2.358}{}$ & $\SI{2.355}{}$ & $\SI{1.886}{}$ & $\SI{1.834}{}$ \\
        $\eta_{\textrm{PtL}}$ & $\SI{}{\percent}$ & $\SI{61.330}{}$ & $\SI{61.844}{}$ & $\SI{57.579}{}$ & $\SI{58.083}{}$ \\
        \midrule
        $\mathit{CAPEX}_{\textrm{sys}}$ & $\SI{}{\EUR\per\year}$ & $\SI{500000}{}$ & $\SI{500000}{}$ & $\SI{500000}{}$ & $\SI{}{500000}$ \\
        $\mathit{CAPEX}_{\textrm{HEN}}$ & $\SI{}{\EUR\per\year}$ & $\SI{39680}{}$ & $\SI{40540}{}$ & $\SI{45980}{}$ & $\SI{25990}{}$ \\
        $\mathit{OPEX}$ & $\SI{}{\EUR\per\year}$ & $\SI{200660}{}$ & $\SI{199140}{}$ & $\SI{260860}{}$ & $\SI{259880}{}$ \\
        $\mathit{TAC}$ & $\SI{}{\EUR\per\year}$ & $\SI{740340}{}$ & $\SI{739680}{}$ & $\SI{806640}{}$ & $\SI{785870}{}$ \\
        \midrule
        $n_{\textrm{HEX}}$ & & $\SI{27}{}$ & $\SI{27}{}$ & $\SI{34}{}$ & $\SI{20}{}$ \\
        $\sum\limits_{i}q_{\mathrm{cu},i}$ & $\SI{}{\kilo\watt}$ & $\SI{222.25}{}$ & $\SI{155.62}{}$ & $\SI{378.63}{}$ & $\SI{255.76}{}$ \\
        $\sum\limits_{j}q_{\mathrm{hu},j}$ & $\SI{}{\kilo\watt}$ & $\SI{5.90}{}$ & $\SI{0}{}$ & $\SI{0}{}$ & $\SI{0}{}$ \\
        $\sum\limits_{i}\sum\limits_{j}\sum\limits_{k}q_{i,j,k}$ & $\SI{}{\kilo\watt}$ & $\SI{717.34}{}$ & $\SI{723.19}{}$ & $\SI{890.26}{}$ & $\SI{890.24}{}$ \\
        \bottomrule
    \end{tabular}
\end{table}

The shape of the Pareto front, whether convex or non-convex, provides information about the nature of the underlying problem and the trade-offs between the different objectives. A convex Pareto front indicates clear dominant solutions where an improvement in one objective necessarily leads to a degradation of the other. A distinct optimal solution usually combines the best values for all objectives in such cases. In contrast, in this case, a non-convex Pareto front indicates no unambiguously dominant solution and different combinations of objective functions represent equivalent alternatives. Looking at the Pareto front from Figure \ref{fig:reults}, we find that the efficiency remains nearly constant in the low-production cost region at about $\SI{58}{\percent}$. For a cell voltage of $U_{\textrm{cell}}^{\textrm{min}} = \SI{1.305}{\volt}$, the production costs vary between \mbox{$\SI{1.834}{\EUR\per\kilogram}$} and \mbox{$\SI{2.057}{\EUR\per\kilogram}$}, while the efficiency varies between $\SI{57.803}{\percent}$ and $\SI{58.371}{\percent}$. In our case, an identical cell voltage means that the solutions differ only in the HEN design and the parameterization of the CS. It can be concluded that, especially in this range, the influence of the design variables has significant effects on the production costs with negligible effects on efficiency. Similarly, in the high production cost range, a significant increase in efficiency can be seen with a minor increase in production cost. From \mbox{$\SI{2.068}{\EUR\per\kilogram}$} to \mbox{$\SI{2.297}{\EUR\per\kilogram}$}, the PtL efficiency increases by 3.963 percentage points to $\SI{61.785}{\percent}$. However, all solutions show different cell voltages in this interval. It can be concluded that the choice of cell voltage has a more significant impact on the objective functions than the choice of HEN and the parameterization of the CS. In the range of highest production costs between \mbox{$\SI{2.297}{\EUR\per\kilogram}$} and \mbox{$\SI{2.355}{\EUR\per\kilogram}$}, four solutions with a cell voltage of $U_{\textrm{cell}}^{\textrm{min}} = \SI{1.275}{\volt}$ are found. The efficiency varies only by about $\SI{0.096}{\percent}$ in the nearly horizontal zone increasing from $\SI{61.785}{\percent}$ . The production costs, on the other hand, differ by up to \mbox{$\SI{5.863}{\ct\per\kilogram}$}.

The results show that in our use case, the cell voltages of the Pareto optimal solutions do not change in the horizontal section at the beginning and end of the Pareto front. In this region, the HEN design and the parameterization of the CS are of crucial importance. However, in the near-vertical part of the Pareto front, the choice of cell voltage significantly affects the objective functions. The HEN design and the parameterization of the CS have only a minor influence. Therefore, the coupled optimization of the operating parameters and the HEN is essential if the impact of different design parameters cannot be estimated a priori.

%% file: 05_conclusion.tex
\section{Conclusion}
\label{sec:conclusion}

In this paper, we present a method that enables coupled optimization of design and operating parameters in the heat exchanger network (HEN) problem. The coupling is realized by modeling an operating point dependent behavior of stream parameters inlet, outlet temperature and flow capacity. All non-linearities are approximated with piecewise linear approximations to keep the problem tractable and to leverage the potential of fast MILP solvers. The transfer to MILP is done highly efficiently using logarithmic coding. The method is applied to a novel \mbox{$\SI{1}{\mega\watt}$} PtL-process. Selecting the cell voltage of the co-SOEC and the combustion system's (CS) parameterization is crucial for the process design. These complex processes cannot be evaluated purposefully based on only one objective. Therefore, multi-criteria optimization of PtL-efficiency and production costs provides a robust decision basis for comprehensive process assessment.

Our results show that coupled optimization leads to better results than the conventional design approach. Higher efficiencies can be achieved with lower production costs at the same time. By modeling the systems’ operational and design parameters, it is possible to optimize the overall process comprehensively and to calculate a Pareto front that reflects the trade-off of the objective functions. The shape of the Pareto front provides valuable insight into the nature of the problem. In this paper, a non-convex Pareto front was observed, suggesting that no uniquely dominant solution and different combinations of objectives represent equivalent alternatives. Efficiency remained approximately constant in the low-production cost region. This indicates that the influence of HEN and CS parameters significantly affects production costs, while the effects on efficiency are negligible. On the other hand, in the high production cost area, a slight increase in production cost resulted in a significant improvement in efficiency. The choice of cell voltage was found to have a more substantial impact on the objective functions than the choice of HEN and the parameterization of the CS. The results underline the relevance of coupled optimization, especially when the effects of the different operational and design parameters cannot be estimated a priori.

With our method, we close a research gap by coupled optimization of operating and design parameters. Simultaneously considering several objective functions enables a comprehensive analysis of the solution space based on Pareto fronts. This allows synergy effects to be exploited and optimal solutions to be identified highly efficiently. Further research and analysis could focus on implementing additional operational or design variables to improve the overall performance and sustainability of the PtL-plant. With this work, we can contribute to the highly efficient and cost-effective production of synthetic fuels, which will promote timely large-scale industrial production and reduce emissions in the long run.

%% file: 06_statements.tex
\section*{Statements and Declarations}
\subsection*{Funding}
This research was funded by the Austrian Research Promotion Agency (FFG) under grant number 884340 and TU Wien Bibliothek through its Open Access Funding Programme.

\subsection*{Competing Interests}
The authors have no relevant financial or non-financial interests to disclose.

\subsection*{Authors Contributions}
The method presented in this paper was developed by David Huber. Testing and evaluation was done by David Huber. The conceptualization of the paper was done by all authors. The first draft was written by David Huber. All authors contributed to the revision of the initial draft. Funding and supervision was done by René Hofmann. All authors read and approved the final manuscript.

%% file: A_appendix.tex
\section{Solver Time}
\label{sec:solvertime}

\begin{figure}[H]
    \centering
    \includegraphics[width=0.7\textwidth]{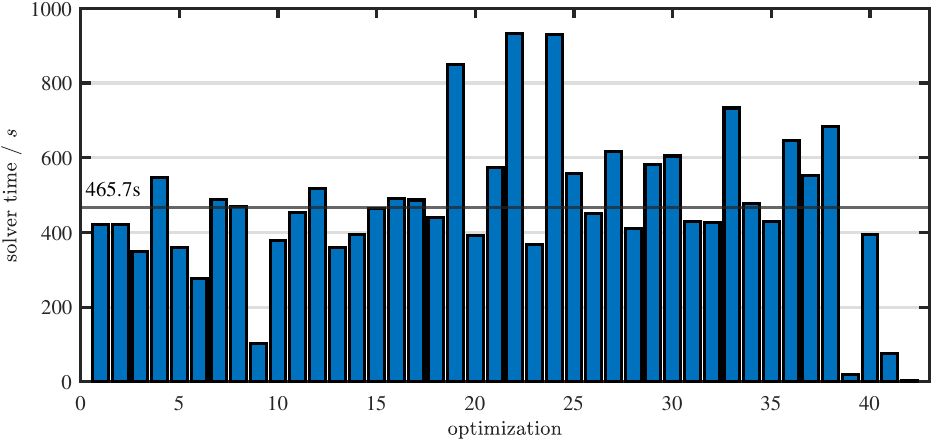}
    \caption{Solver time for each optimization (blue bars) and mean solver time (black line).} \label{fig:solverTime}
\end{figure}

Figure \ref{fig:solverTime} shows the solver time for the $42$ points of the Pareto front as bar plot. The solver time refers to the time it takes the solver to find a solution within the defined MIP gap of $\SI{1}{\percent}$. The mean solver time was $\SI{465.7}{\second}$. The total solver time was $8.5809 \cdot 10^4 \SI{}{\second}$ or $\SI{23.8}{\hour}$.

\section{Stream Plots}
\label{sec:streamplots}
\begin{figure}[H]
    \centering
    \includegraphics[width=1\textwidth]{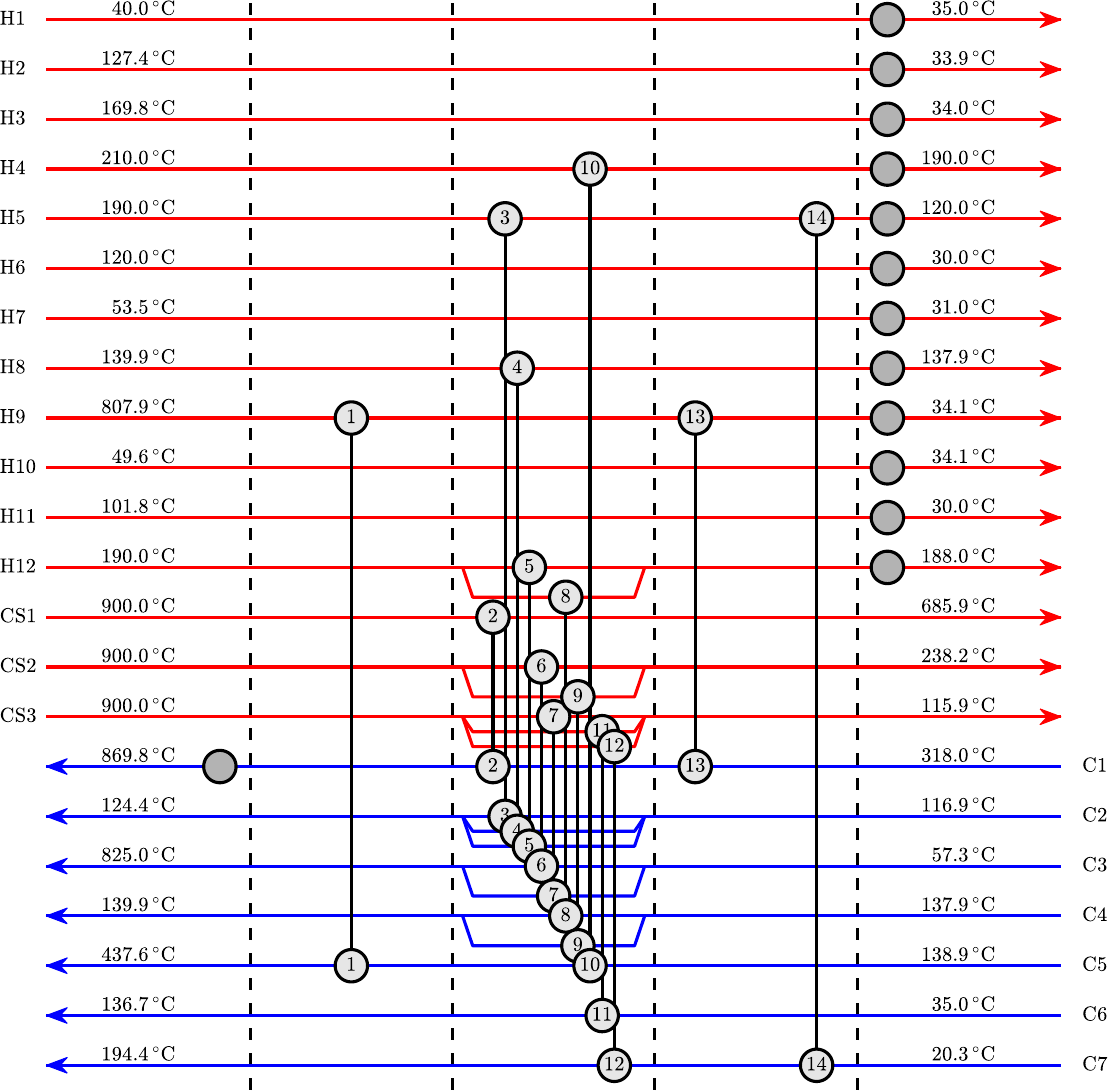}
    \caption{Stream plot with $27$ heat exchangers resulting from the empirical design without optimization. Characteristic figures: $\eta_{\textrm{PtL}} = \SI{61.330}{\percent}$, $c_{\textrm{prod}} = \SI{2.358}{\EUR\per\kilogram}$, $U_{\textrm{cell}} = \SI{1.275}{\volt}$.} \label{fig:streamplotempirical1275}
\end{figure}

\begin{figure}[H]
    \centering
    \includegraphics[width=1\textwidth]{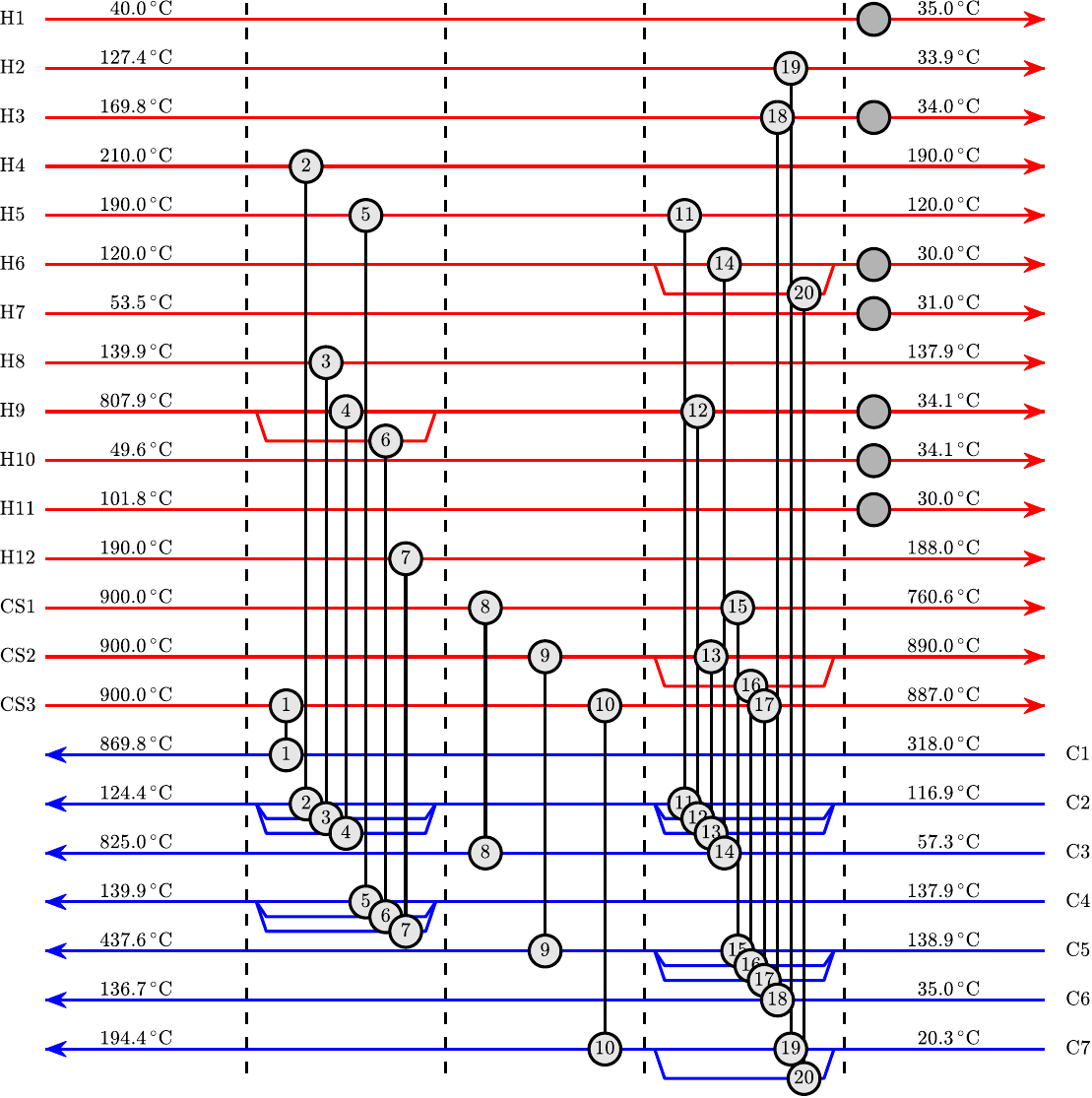}
    \caption{Stream plot with $27$ heat exchangers resulting from the coupled optimization. Characteristic figures: $\eta_{\textrm{PtL}} = \SI{61.844}{\percent}$, $c_{\textrm{prod}} = \SI{2.355}{\EUR\per\kilogram}$, $U_{\textrm{cell}} = \SI{1.275}{\volt}$.} \label{fig:streamplotopt1275}
\end{figure}

\begin{figure}[H]
    \centering
    \includegraphics[width=1\textwidth]{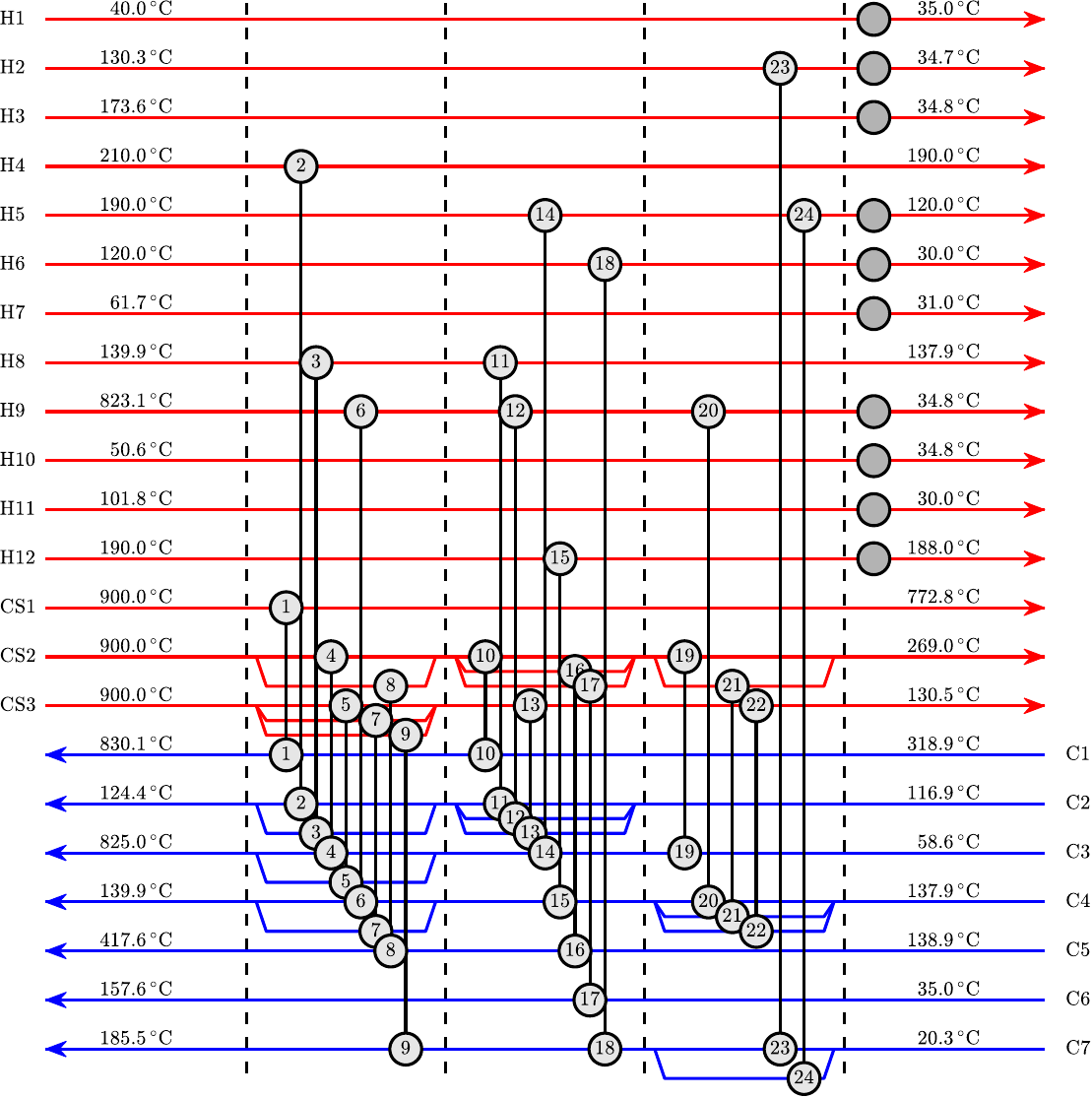}
    \caption{Stream plot with $34$ heat exchangers resulting from the empirical design without optimization. Characteristic figures: $\eta_{\textrm{PtL}} = \SI{57.579}{\percent}$, $c_{\textrm{prod}} = \SI{1.886}{\EUR\per\kilogram}$, $U_{\textrm{cell}} = \SI{1.305}{\volt}$.} \label{fig:streamplotempirical1305}
\end{figure}

\begin{figure}[H]
    \centering
    \includegraphics[width=1\textwidth]{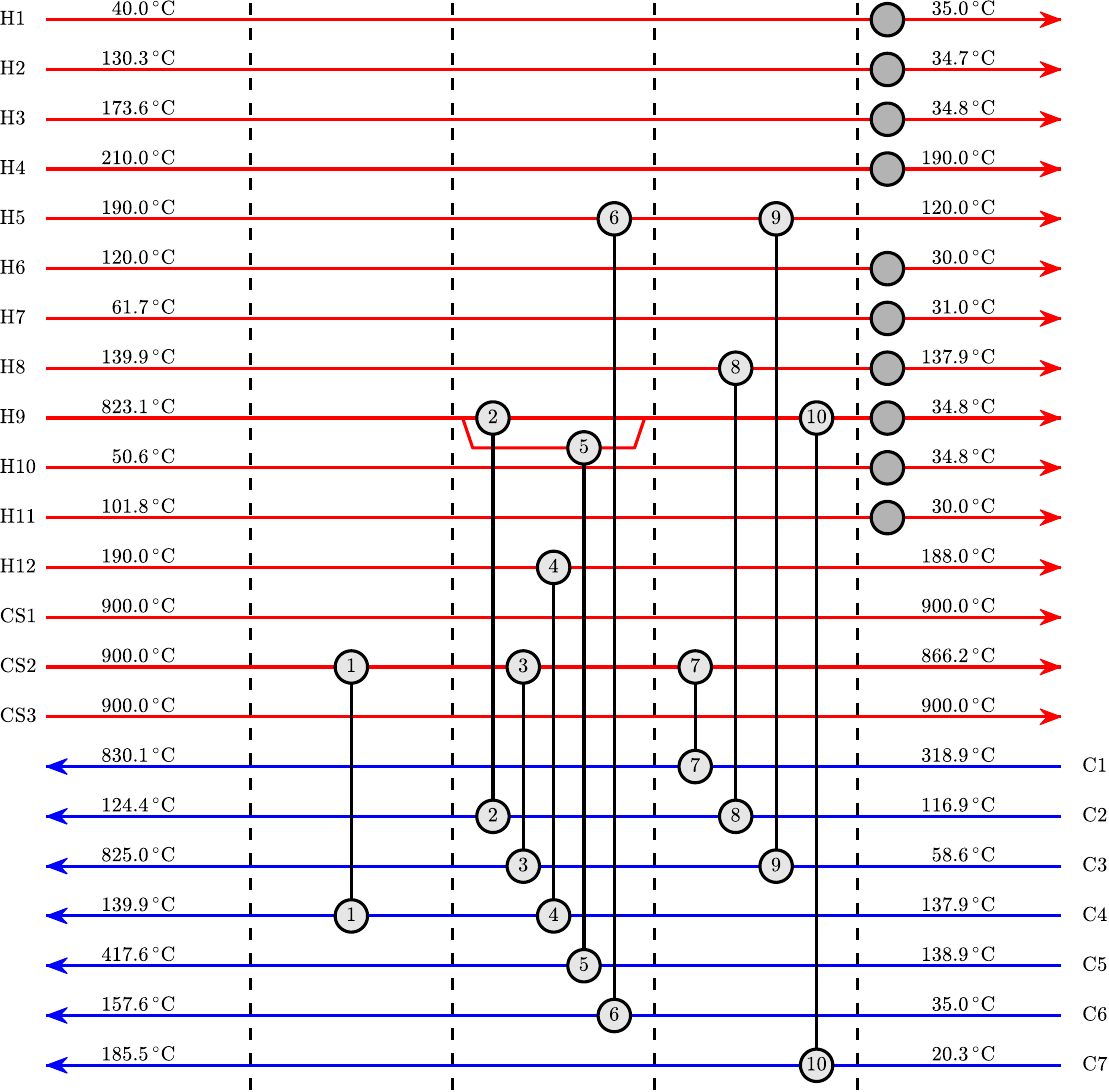}
    \caption{Stream plot with $20$ heat exchangers resulting from the empirical design without optimization. Characteristic figures: $\eta_{\textrm{PtL}} = \SI{58.083}{\percent}$, $c_{\textrm{prod}} = \SI{1.834}{\EUR\per\kilogram}$, $U_{\textrm{cell}} = \SI{1.305}{\volt}$.} \label{fig:streamplotopt13005}
\end{figure}

Figures \ref{fig:streamplotempirical1275} and \ref{fig:streamplotempirical1305} show the empirically constructed HENs for a cell voltage of $U_{\textrm{cell}} = \SI{1.275}{\volt}$ and $U_{\textrm{cell}} = \SI{1.305}{\volt}$, respectively. In addition, Figures \ref{fig:streamplotopt1275} and \ref{fig:streamplotopt13005} show the HENs with the same cell voltage resulting from the coupled optimization. Based on the detailed results from Table \ref{tab:results}, the corner points of the Pareto front have been used.

With the coupled optimization, the design of the combustion system (CS) is coupled with the optimization of the HEN, among other things. Accordingly, the outlet temperatures and flow capacities of the corresponding streams CS1-3 are implemented as variables. It can be seen that compared to the empirically defined designs, the outlet temperatures of stream CS1-3 are significantly higher. Du to the higher temperature differences, the heat exchanger surfaces can be reduced at the same performance, lowering costs.

%% file: nomenclature.tex
\makenomenclature

\renewcommand\nomgroup[1]{%
  \item[\bfseries
  \ifstrequal{#1}{A}{Acronyms}{%
  \ifstrequal{#1}{V}{Variables}{%
  \ifstrequal{#1}{L}{Subscripts}{%
  \ifstrequal{#1}{H}{Superscripts}{%
  \ifstrequal{#1}{S}{Sets}}}}}%
]}

\setlength{\nomlabelwidth}{1.5cm}

\newcommand{\nomunit}[1]{%
\renewcommand{\nomentryend}{\hspace*{\fill}#1}}

\nomenclature[A]{SOEC}{solid oxide electrolysis cell}

\nomenclature[A]{MILP}{mixed-integer linear programming}

\nomenclature[A]{PtL}{power-to-liquid}

\nomenclature[A]{HEN}{heat exchanger network}

\nomenclature[A]{HENS}{heat exchanger network synthesis}

\nomenclature[A]{DAC}{direct air capture}

\nomenclature[A]{FT}{fischer-tropsch}

\nomenclature[A]{PEM}{polymer electrolyte membrane}

\nomenclature[A]{GA}{genetic algorithm}

\nomenclature[A]{LCOP}{levelized cost of product}

\nomenclature[A]{IFE}{Innovation Flüssige Energie, eng.: Innovation Liquid Energy}

\nomenclature[A]{CS}{combustion system}

\nomenclature[A]{RMSE}{root-mean-square error}

\nomenclature[A]{CAPEX}{annual capital expenses}

\nomenclature[A]{OPEX}{operational expenditures}

\nomenclature[A]{LMTD}{logarithmic mean temperature difference}

\nomenclature[A]{HEX}{heat exchanger}

\nomenclature[A]{PV}{photovoltaics}

\nomenclature[V]{$F$}{{flow capacity}\nomunit{\SI{}{\kilo\watt\per\kelvin}}}

\nomenclature[V]{$U_{\textrm{cell}}$}{{cell voltage}
\nomunit{\SI{}{\volt}}}

\nomenclature[V]{$\dot{m}$}{{mass flow}
\nomunit{\SI{}{\kg\per\hour}}}

\nomenclature[V]{$T$}{{temperature}
\nomunit{\SI{}{\celsius}}}

\nomenclature[V]{$h_{\textrm{prod}}$}{{specific enthalpy of the product}
\nomunit{\SI{}{\mega\joule\per\kilogram}}}

\nomenclature[V]{$\rho$}{{density}
\nomunit{\SI{}{\kilogram\cubic\meter}}}

\nomenclature[V]{$\mu$}{{dynamic viscosity}
\nomunit{\SI{}{\milli\pascal\per\second}}}

\nomenclature[V]{$N_{\textrm{st}}$}{number of stages of the HEN superstructure}

\nomenclature[V]{$n_{\textrm{HEX}}$}{number of HEX in the HEN superstructure}

\nomenclature[V]{$P_{\textrm{sys}}$}{{electrical energy demand w/o utilities}
\nomunit{\SI{}{\kW}}}

\nomenclature[V]{$U$}{{overall heat transfer coefficient}\nomunit{\SI{}{\kilo\watt\per\meter\squared\per\kelvin}}}

\nomenclature[V]{$\eta_{\textrm{PtL}}$}{{PtL-efficiency}
\nomunit{\SI{}{\percent}}}

\nomenclature[V]{$P_{\textrm{el}}$}{{total electrical energy demand}
\nomunit{\SI{}{\kilo\watt}}}

\nomenclature[V]{$\dot{H}$}{{chemically bounded energy in FT-products}
\nomunit{\SI{}{\kilo\watt}}}

\nomenclature[V]{$\varepsilon$}{{coefficient of performance}
\nomunit{\SI{}{}}}

\nomenclature[V]{$q$}{{heat flow}
\nomunit{\SI{}{\kilo\watt}}}

\nomenclature[V]{$t$}{{annual full load hours}
\nomunit{\SI{}{\hour\per\year}}}

\nomenclature[V]{$C_{\textrm{sys}}$}{{investment costs}
\nomunit{\SI{}{\EUR}}}

\nomenclature[V]{$a$}{{depreciation period}
\nomunit{\SI{}{\year}}}

\nomenclature[V]{$AF_{\textrm{inv}}$}{{investment annualization factor}
\nomunit{\SI{}{\per\year}}}

\nomenclature[V]{$AF_{\textrm{op}}$}{{operational annualization factor}
\nomunit{\SI{}{}}}

\nomenclature[V]{$c_{\textrm{prod}}$}{{product costs}
\nomunit{\SI{}{\EUR\per\kilogram}}}

\nomenclature[V]{$c_{\text{f}}$}{{feedstock costs}
\nomunit{\SI{}{\EUR\per\tonne}}}

\nomenclature[V]{$c_{\textrm{el}}$}{{electricity costs}
\nomunit{\SI{}{\EUR\per\mega\watt\per\hour}}}

\nomenclature[V]{$c_{\textrm{f,hex}}$}{{step-fixed HEX costs}
\nomunit{\SI{}{\EUR\per\year}}}

\nomenclature[V]{$c_{\textrm{v,hex}}$}{{variable HEX costs}
\nomunit{\SI{}{\EUR\per\betaCo\per\year}}}

\nomenclature[V]{$\beta$}{{cost exponent}\nomunit{\SI{}{}}}

\nomenclature[V]{$z$}{{binary variable for existence of HEX}\nomunit{\SI{}{}}}

\nomenclature[V]{$\Delta T_{\textrm{min}}$}{{minimum temperature difference}
\nomunit{\SI{}{\kelvin}}}

\nomenclature[V]{$\mathit{TAC}$}{{total annual costs}
\nomunit{\SI{}{\EUR\per\year}}}

\nomenclature[V]{$\mathit{LMTD}$}{{logarithmic mean temperature difference}
\nomunit{\SI{}{\kelvin}}}

\nomenclature[L]{$\textrm{prod}$}{production}

\nomenclature[L]{$\textrm{cu}$}{cold utility}

\nomenclature[L]{$\textrm{hu}$}{hot utility}

\nomenclature[L]{$\textrm{hex}$}{heat exchanger}

\nomenclature[L]{$\textrm{el}$}{electric}

\nomenclature[L]{$i$}{hot stream}

\nomenclature[L]{$j$}{cold stream}

\nomenclature[L]{$k$}{stage of the HEN superstructure}

\nomenclature[H]{in}{inlet}

\nomenclature[H]{out}{outlet}

\nomenclature[H]{min}{minimum}

\nomenclature[H]{max}{maximum}

\printnomenclature

%% file: _main_paper.bbl
\begin{thebibliography}{37}
\expandafter\ifx\csname natexlab\endcsname\relax\def\natexlab#1{#1}\fi
\providecommand{\url}[1]{\texttt{#1}}
\providecommand{\href}[2]{#2}
\providecommand{\path}[1]{#1}
\providecommand{\DOIprefix}{doi:}
\providecommand{\ArXivprefix}{arXiv:}
\providecommand{\URLprefix}{URL: }
\providecommand{\Pubmedprefix}{pmid:}
\providecommand{\doi}[1]{\href{http://dx.doi.org/#1}{\path{#1}}}
\providecommand{\Pubmed}[1]{\href{pmid:#1}{\path{#1}}}
\providecommand{\bibinfo}[2]{#2}
\ifx\xfnm\relax \def\xfnm[#1]{\unskip,\space#1}\fi
\bibitem[{{International Energy Agency}(2022{\natexlab{a}})}]{international_energy_agency_co2_2022}
\bibinfo{author}{{International Energy Agency}}, \bibinfo{title}{{CO2} {Emissions} in 2022: {Global} {Energy} {Review}}, \bibinfo{type}{Technical Report}, \bibinfo{year}{2022}{\natexlab{a}}. \URLprefix \url{https://www.iea.org/reports/co2-emissions-in-2022}.
\bibitem[{{International Energy Agency}(2022{\natexlab{b}})}]{international_energy_agency_global_2022}
\bibinfo{author}{{International Energy Agency}}, \bibinfo{title}{Global {EV} {Outlook} 2022}, \bibinfo{type}{Technical Report}, \bibinfo{year}{2022}{\natexlab{b}}. \URLprefix \url{https://www.iea.org/reports/global-ev-outlook-2022}.
\bibitem[{Hänggi et~al.(2019)Hänggi, Elbert, Bütler, Cabalzar, Teske, Bach, and Onder}]{hanggi_review_2019}
\bibinfo{author}{S.~Hänggi}, \bibinfo{author}{P.~Elbert}, \bibinfo{author}{T.~Bütler}, \bibinfo{author}{U.~Cabalzar}, \bibinfo{author}{S.~Teske}, \bibinfo{author}{C.~Bach}, \bibinfo{author}{C.~Onder},
\newblock \bibinfo{title}{A review of synthetic fuels for passenger vehicles},
\newblock \bibinfo{journal}{Energy Reports} \bibinfo{volume}{5} (\bibinfo{year}{2019}) \bibinfo{pages}{555--569}. \DOIprefix\doi{10.1016/j.egyr.2019.04.007}.
\bibitem[{Pinsky et~al.(2020)Pinsky, Sabharwall, Hartvigsen, and O’Brien}]{pinsky_comparative_2020}
\bibinfo{author}{R.~Pinsky}, \bibinfo{author}{P.~Sabharwall}, \bibinfo{author}{J.~Hartvigsen}, \bibinfo{author}{J.~O’Brien},
\newblock \bibinfo{title}{Comparative review of hydrogen production technologies for nuclear hybrid energy systems},
\newblock \bibinfo{journal}{Progress in Nuclear Energy} \bibinfo{volume}{123} (\bibinfo{year}{2020}) \bibinfo{pages}{103317}. \DOIprefix\doi{10.1016/j.pnucene.2020.103317}.
\bibitem[{Marchese et~al.(2020)Marchese, Giglio, Santarelli, and Lanzini}]{marchese_energy_2020}
\bibinfo{author}{M.~Marchese}, \bibinfo{author}{E.~Giglio}, \bibinfo{author}{M.~Santarelli}, \bibinfo{author}{A.~Lanzini},
\newblock \bibinfo{title}{Energy performance of {Power}-to-{Liquid} applications integrating biogas upgrading, reverse water gas shift, solid oxide electrolysis and {Fischer}-{Tropsch} technologies},
\newblock \bibinfo{journal}{Energy Conversion and Management: X} \bibinfo{volume}{6} (\bibinfo{year}{2020}) \bibinfo{pages}{100041}. \DOIprefix\doi{10.1016/j.ecmx.2020.100041}.
\bibitem[{{Siemens energy}(2023)}]{siemens_energy_haru_2023}
\bibinfo{author}{{Siemens energy}}, \bibinfo{title}{Haru {Oni}: {Base} camp of the future}, \bibinfo{year}{2023}. \URLprefix \url{https://www.siemens-energy.com/global/en/news/magazine/2022/haru-oni.html}.
\bibitem[{{Norsk e-fuel}(2023)}]{norsk_e-fuel_driving_2023}
\bibinfo{author}{{Norsk e-fuel}}, \bibinfo{title}{Driving the transition to renewable aviation today}, \bibinfo{year}{2023}. \URLprefix \url{https://www.norsk-e-fuel.com/}.
\bibitem[{Marlin et~al.(2018)Marlin, Sarron, and Sigurbjörnsson}]{marlin_process_2018}
\bibinfo{author}{D.~S. Marlin}, \bibinfo{author}{E.~Sarron}, \bibinfo{author}{O.~Sigurbjörnsson},
\newblock \bibinfo{title}{Process {Advantages} of {Direct} {CO2} to {Methanol} {Synthesis}},
\newblock \bibinfo{journal}{Frontiers in Chemistry} \bibinfo{volume}{6} (\bibinfo{year}{2018}) \bibinfo{pages}{446}. \DOIprefix\doi{10.3389/fchem.2018.00446}.
\bibitem[{{INERATEC}(2022)}]{ineratec_industrial_2022}
\bibinfo{author}{{INERATEC}}, \bibinfo{title}{Industrial {Power}-to-{Liquid} {Pioneer} {Plant} in {Germany}}, \bibinfo{year}{2022}. \URLprefix \url{https://planet-a.com/startups/ineratec/}.
\bibitem[{Zhao et~al.(2022)Zhao, Yu, Ren, Makowski, Granat, Nahorski, and Ma}]{zhao_how_2022}
\bibinfo{author}{J.~Zhao}, \bibinfo{author}{Y.~Yu}, \bibinfo{author}{H.~Ren}, \bibinfo{author}{M.~Makowski}, \bibinfo{author}{J.~Granat}, \bibinfo{author}{Z.~Nahorski}, \bibinfo{author}{T.~Ma},
\newblock \bibinfo{title}{How the power-to-liquid technology can contribute to reaching carbon neutrality of the {China}'s transportation sector?},
\newblock \bibinfo{journal}{Energy} \bibinfo{volume}{261} (\bibinfo{year}{2022}) \bibinfo{pages}{125058}. \DOIprefix\doi{10.1016/j.energy.2022.125058}.
\bibitem[{Ngan Do et~al.(2022)Ngan Do, You, and Kim}]{ngando_co_2022}
\bibinfo{author}{T.~Ngan Do}, \bibinfo{author}{C.~You}, \bibinfo{author}{J.~Kim},
\newblock \bibinfo{title}{A {CO} 2 utilization framework for liquid fuels and chemical production: techno-economic and environmental analysis},
\newblock \bibinfo{journal}{Energy \& Environmental Science} \bibinfo{volume}{15} (\bibinfo{year}{2022}) \bibinfo{pages}{169--184}. \DOIprefix\doi{10.1039/D1EE01444G}, \bibinfo{note}{publisher: Royal Society of Chemistry}.
\bibitem[{Ueckerdt et~al.(2021)Ueckerdt, Bauer, Dirnaichner, Everall, Sacchi, and Luderer}]{ueckerdt_potential_2021}
\bibinfo{author}{F.~Ueckerdt}, \bibinfo{author}{C.~Bauer}, \bibinfo{author}{A.~Dirnaichner}, \bibinfo{author}{J.~Everall}, \bibinfo{author}{R.~Sacchi}, \bibinfo{author}{G.~Luderer},
\newblock \bibinfo{title}{Potential and risks of hydrogen-based e-fuels in climate change mitigation},
\newblock \bibinfo{journal}{Nature Climate Change} \bibinfo{volume}{11} (\bibinfo{year}{2021}) \bibinfo{pages}{384--393}. \DOIprefix\doi{10.1038/s41558-021-01032-7}, \bibinfo{note}{number: 5 Publisher: Nature Publishing Group}.
\bibitem[{Ram et~al.(2020)Ram, Galimova, Bogdanov, Fasihi, Gulagi, Micheli, Crone, and Breyer}]{ram_powerfuels_2020}
\bibinfo{author}{M.~Ram}, \bibinfo{author}{T.~Galimova}, \bibinfo{author}{D.~Bogdanov}, \bibinfo{author}{M.~Fasihi}, \bibinfo{author}{A.~Gulagi}, \bibinfo{author}{M.~Micheli}, \bibinfo{author}{K.~Crone}, \bibinfo{author}{C.~Breyer},
\newblock \bibinfo{title}{Powerfuels in a {Renewable} {Energy} {World}: {Global} {Volumes}, {Costs}, and {Trading} 2030 to 2050}  (\bibinfo{year}{2020}). \DOIprefix\doi{10.13140/RG.2.2.32687.56487}, \bibinfo{note}{publisher: Lappeenranta – Lahti University of Technology LUT}.
\bibitem[{{U.S. Energy Information Administration}(2023)}]{us_energy_information_administration_us_2023}
\bibinfo{author}{{U.S. Energy Information Administration}}, \bibinfo{title}{U.{S}. {Total} {Gasoline} {Wholesale}/{Resale} {Price} by {Refiners} ({Dollars} per {Gallon})}, \bibinfo{year}{2023}. \URLprefix \url{https://www.eia.gov/petroleum/data.php}.
\bibitem[{Ueckerdt and Odenweller(2023)}]{ueckerdt_e-fuels_2023}
\bibinfo{author}{F.~Ueckerdt}, \bibinfo{author}{A.~Odenweller}, \bibinfo{title}{E-{Fuels} - {Aktueller} {Stand} und {Projektionen}}, \bibinfo{type}{Technical Report}, Potsdam-Institut für Klimafolgenforschung, \bibinfo{year}{2023}.
\bibitem[{Al-Rashed and Afrand(2021)}]{al-rashed_multi-criteria_2021}
\bibinfo{author}{A.~A. Al-Rashed}, \bibinfo{author}{M.~Afrand},
\newblock \bibinfo{title}{Multi-criteria exergoeconomic optimization for a combined gas turbine-supercritical {CO2} plant with compressor intake cooling fueled by biogas from anaerobic digestion},
\newblock \bibinfo{journal}{Energy} \bibinfo{volume}{223} (\bibinfo{year}{2021}) \bibinfo{pages}{119997}. \DOIprefix\doi{10.1016/j.energy.2021.119997}.
\bibitem[{Cao et~al.(2022)Cao, Dhahad, Hussen, Anqi, Farouk, and Issakhov}]{cao_development_2022}
\bibinfo{author}{Y.~Cao}, \bibinfo{author}{H.~A. Dhahad}, \bibinfo{author}{H.~M. Hussen}, \bibinfo{author}{A.~E. Anqi}, \bibinfo{author}{N.~Farouk}, \bibinfo{author}{A.~Issakhov},
\newblock \bibinfo{title}{Development and tri-objective optimization of a novel biomass to power and hydrogen plant: {A} comparison of fueling with biomass gasification or biomass digestion},
\newblock \bibinfo{journal}{Energy} \bibinfo{volume}{238} (\bibinfo{year}{2022}) \bibinfo{pages}{122010}. \DOIprefix\doi{10.1016/j.energy.2021.122010}.
\bibitem[{Hai et~al.(2023)Hai, Hikmat Hama~Aziz, Zhou, Dhahad, Sharma, Fahad~Almojil, Ibrahim~Almohana, Fahmi~Alali, Ismail~Kh, Mehrez, and Abdelrahman}]{hai_-neural_2023}
\bibinfo{author}{T.~Hai}, \bibinfo{author}{K.~Hikmat Hama~Aziz}, \bibinfo{author}{J.~Zhou}, \bibinfo{author}{H.~A. Dhahad}, \bibinfo{author}{K.~Sharma}, \bibinfo{author}{S.~Fahad~Almojil}, \bibinfo{author}{A.~Ibrahim~Almohana}, \bibinfo{author}{A.~Fahmi~Alali}, \bibinfo{author}{T.~Ismail~Kh}, \bibinfo{author}{S.~Mehrez}, \bibinfo{author}{A.~Abdelrahman},
\newblock \bibinfo{title}{-{Neural} network-based optimization of hydrogen fuel production energy system with proton exchange electrolyzer supported nanomaterial},
\newblock \bibinfo{journal}{Fuel} \bibinfo{volume}{332} (\bibinfo{year}{2023}) \bibinfo{pages}{125827}. \DOIprefix\doi{10.1016/j.fuel.2022.125827}.
\bibitem[{Wang et~al.(2023)Wang, Guo, Pei, He, Liu, and Li}]{wang_rolling_2023}
\bibinfo{author}{C.~Wang}, \bibinfo{author}{S.~Guo}, \bibinfo{author}{H.~Pei}, \bibinfo{author}{Y.~He}, \bibinfo{author}{D.~Liu}, \bibinfo{author}{M.~Li},
\newblock \bibinfo{title}{Rolling optimization based on holism for the operation strategy of solar tower power plant},
\newblock \bibinfo{journal}{Applied Energy} \bibinfo{volume}{331} (\bibinfo{year}{2023}) \bibinfo{pages}{120473}. \DOIprefix\doi{10.1016/j.apenergy.2022.120473}.
\bibitem[{Linnhoff(1998)}]{linnhoff_introduction_1998}
\bibinfo{author}{M.~Linnhoff},
\newblock \bibinfo{title}{Introduction to pinch technology},
\newblock \bibinfo{journal}{Targeting House, Gadbrook Park, Northwich, Cheshire, CW9 7UZ, England}  (\bibinfo{year}{1998}).
\bibitem[{Yee and Grossmann(1990)}]{yee_simultaneous_1990}
\bibinfo{author}{T.~F. Yee}, \bibinfo{author}{I.~E. Grossmann},
\newblock \bibinfo{title}{Simultaneous optimization models for heat integration—{II}. {Heat} exchanger network synthesis},
\newblock \bibinfo{journal}{Computers \& Chemical Engineering} \bibinfo{volume}{14} (\bibinfo{year}{1990}) \bibinfo{pages}{1165--1184}. \DOIprefix\doi{10.1016/0098-1354(90)85010-8}.
\bibitem[{Escobar and Trierweiler(2013)}]{escobar_optimal_2013}
\bibinfo{author}{M.~Escobar}, \bibinfo{author}{J.~O. Trierweiler},
\newblock \bibinfo{title}{Optimal heat exchanger network synthesis: {A} case study comparison},
\newblock \bibinfo{journal}{Applied Thermal Engineering} \bibinfo{volume}{51} (\bibinfo{year}{2013}) \bibinfo{pages}{801--826}. \DOIprefix\doi{10.1016/j.applthermaleng.2012.10.022}.
\bibitem[{Liu et~al.(2022)Liu, Yang, Yang, and Qian}]{liu_extended_2022}
\bibinfo{author}{Z.~Liu}, \bibinfo{author}{L.~Yang}, \bibinfo{author}{S.~Yang}, \bibinfo{author}{Y.~Qian},
\newblock \bibinfo{title}{An extended stage-wise superstructure for heat exchanger network synthesis with intermediate placement of multiple utilities},
\newblock \bibinfo{journal}{Energy} \bibinfo{volume}{248} (\bibinfo{year}{2022}) \bibinfo{pages}{123372}. \DOIprefix\doi{10.1016/j.energy.2022.123372}.
\bibitem[{Huber et~al.(2023)Huber, Birkelbach, and Hofmann}]{huber_hens_2023}
\bibinfo{author}{D.~Huber}, \bibinfo{author}{F.~Birkelbach}, \bibinfo{author}{R.~Hofmann},
\newblock \bibinfo{title}{{HENS} {Unchained}: {MILP} {Implementation} of {Multi}-{Stage} {Utilities} with {Stream} {Splits}, {Variable} {Temperatures} and {Flow} {Capacities}},
\newblock \bibinfo{journal}{Energies} \bibinfo{volume}{16} (\bibinfo{year}{2023}) \bibinfo{pages}{4732}. \DOIprefix\doi{10.3390/en16124732}.
\bibitem[{Vielma and Nemhauser(2011)}]{vielma_modeling_2011}
\bibinfo{author}{J.~P. Vielma}, \bibinfo{author}{G.~L. Nemhauser},
\newblock \bibinfo{title}{Modeling disjunctive constraints with a logarithmic number of binary variables and constraints},
\newblock \bibinfo{journal}{Mathematical Programming} \bibinfo{volume}{128} (\bibinfo{year}{2011}) \bibinfo{pages}{49--72}. \DOIprefix\doi{10.1007/s10107-009-0295-4}.
\bibitem[{Halmschlager et~al.(2022)Halmschlager, Beck, Knöttner, Koller, and Hofmann}]{Halmschlager2022}
\bibinfo{author}{D.~Halmschlager}, \bibinfo{author}{A.~Beck}, \bibinfo{author}{S.~Knöttner}, \bibinfo{author}{M.~Koller}, \bibinfo{author}{R.~Hofmann},
\newblock \bibinfo{title}{Combined optimization for retrofitting of heat recovery and thermal energy supply in industrial systems},
\newblock \bibinfo{journal}{Applied Energy} \bibinfo{volume}{305} (\bibinfo{year}{2022}) \bibinfo{pages}{117820}. \DOIprefix\doi{10.1016/j.apenergy.2021.117820}.
\bibitem[{Adelung(2022)}]{adelung_global_2022}
\bibinfo{author}{S.~Adelung},
\newblock \bibinfo{title}{Global sensitivity and uncertainty analysis of a {Fischer}-{Tropsch} based {Power}-to-{Liquid} process},
\newblock \bibinfo{journal}{Journal of CO2 Utilization} \bibinfo{volume}{65} (\bibinfo{year}{2022}) \bibinfo{pages}{102171}. \DOIprefix\doi{10.1016/j.jcou.2022.102171}.
\bibitem[{Herz et~al.(2021)Herz, Rix, Jacobasch, Müller, Reichelt, Jahn, and Michaelis}]{herz_economic_2021}
\bibinfo{author}{G.~Herz}, \bibinfo{author}{C.~Rix}, \bibinfo{author}{E.~Jacobasch}, \bibinfo{author}{N.~Müller}, \bibinfo{author}{E.~Reichelt}, \bibinfo{author}{M.~Jahn}, \bibinfo{author}{A.~Michaelis},
\newblock \bibinfo{title}{Economic assessment of {Power}-to-{Liquid} processes – {Influence} of electrolysis technology and operating conditions},
\newblock \bibinfo{journal}{Applied Energy} \bibinfo{volume}{292} (\bibinfo{year}{2021}) \bibinfo{pages}{116655}. \DOIprefix\doi{10.1016/j.apenergy.2021.116655}.
\bibitem[{Dieterich et~al.(2020)Dieterich, Buttler, Hanel, Spliethoff, and Fendt}]{dieterich_power--liquid_2020}
\bibinfo{author}{V.~Dieterich}, \bibinfo{author}{A.~Buttler}, \bibinfo{author}{A.~Hanel}, \bibinfo{author}{H.~Spliethoff}, \bibinfo{author}{S.~Fendt},
\newblock \bibinfo{title}{Power-to-liquid \textit{via} synthesis of methanol, {DME} or {Fischer}–{Tropsch}-fuels: a review},
\newblock \bibinfo{journal}{Energy \& Environmental Science} \bibinfo{volume}{13} (\bibinfo{year}{2020}) \bibinfo{pages}{3207--3252}. \DOIprefix\doi{10.1039/D0EE01187H}.
\bibitem[{König et~al.(2015)König, Freiberg, Dietrich, and Wörner}]{konig_techno-economic_2015}
\bibinfo{author}{D.~H. König}, \bibinfo{author}{M.~Freiberg}, \bibinfo{author}{R.-U. Dietrich}, \bibinfo{author}{A.~Wörner},
\newblock \bibinfo{title}{Techno-economic study of the storage of fluctuating renewable energy in liquid hydrocarbons},
\newblock \bibinfo{journal}{Fuel} \bibinfo{volume}{159} (\bibinfo{year}{2015}) \bibinfo{pages}{289--297}. \DOIprefix\doi{10.1016/j.fuel.2015.06.085}.
\bibitem[{{EurEau}(2021)}]{eureau_europes_2021}
\bibinfo{author}{{EurEau}}, \bibinfo{title}{Europe’s {Water} in {Figures}: {An} overview of the {European} drinking water and waste water sectors}, \bibinfo{year}{2021}. \URLprefix \url{https://www.eureau.org/resources/publications/eureau-publications/5824-europe-s-water-in-figures-2021/file}.
\bibitem[{IETA(2022)}]{ieta_ghg_2022}
\bibinfo{author}{IETA}, \bibinfo{title}{{GHG} {Market} {Sentiment} {Survey} 2022}, \bibinfo{type}{Technical Report}, \bibinfo{year}{2022}. \URLprefix \url{https://www.ieta.org/resources/Documents/IETA%20GHG%20Market%20Sentiment%20Survey%20Report%202022.pdf}.
\bibitem[{Janssen et~al.(2022)Janssen, Weeda, Detz, and van~der Zwaan}]{janssen_country-specific_2022}
\bibinfo{author}{J.~Janssen}, \bibinfo{author}{M.~Weeda}, \bibinfo{author}{R.~J. Detz}, \bibinfo{author}{B.~van~der Zwaan},
\newblock \bibinfo{title}{Country-specific cost projections for renewable hydrogen production through off-grid electricity systems},
\newblock \bibinfo{journal}{Applied Energy} \bibinfo{volume}{309} (\bibinfo{year}{2022}) \bibinfo{pages}{118398}. \DOIprefix\doi{10.1016/j.apenergy.2021.118398}.
\bibitem[{noa(2021)}]{noauthor_dace_2021}
\bibinfo{title}{{DACE} price booklet: cost information for estimation and comparison}, \bibinfo{edition}{edition 35} ed., \bibinfo{publisher}{DACE Cost and Value}, \bibinfo{address}{Nijkerk}, \bibinfo{year}{2021}. \bibinfo{note}{OCLC: 1303572720}.
\bibitem[{{The MathWorks Inc.}(2022)}]{the_mathworks_inc_matlab_2022}
\bibinfo{author}{{The MathWorks Inc.}}, \bibinfo{title}{{MATLAB} version: 9.13.0 ({R2022b})}, \bibinfo{year}{2022}. \URLprefix \url{https://www.mathworks.com}.
\bibitem[{Löfberg(2004)}]{lofberg_toolbox_2004}
\bibinfo{author}{J.~Löfberg},
\newblock \bibinfo{title}{A toolbox for modeling and optimization in {MATLAB}},
\newblock \bibinfo{journal}{Proceedings of the CACSD Conference}  (\bibinfo{year}{2004}) \bibinfo{pages}{289}. \DOIprefix\doi{10.1109/CACSD.2004.1393890}.
\bibitem[{{Gurobi Optimization, LLC}(2023)}]{gurobi_optimization_llc_gurobi_2023}
\bibinfo{author}{{Gurobi Optimization, LLC}}, \bibinfo{title}{Gurobi {Optimizer} {Reference} {Manual}}, \bibinfo{year}{2023}. \URLprefix \url{https://www.gurobi.com}.

\end{thebibliography}
